\documentclass[aps,twocolumn,
superscriptaddress,
footinbib,
prl]{revtex4-1}

\usepackage{amsmath,amssymb,bm}
\usepackage{graphicx}
\usepackage{epstopdf}
\usepackage{latexsym}
\usepackage[normalem]{ulem} 
\usepackage[caption=false]{subfig}
\usepackage{subfig}
\usepackage[usenames,dvipsnames]{color}
\usepackage{hyperref}
\usepackage{natbib}
\usepackage{xcolor}
\usepackage{multirow}
\usepackage{dsfont}
\usepackage{relsize}
\usepackage{float}
\usepackage{algorithm}
\usepackage[noend]{algpseudocode}
\usepackage{listings}

\usepackage{framed}

\usepackage[utf8]{inputenc} 

\begin{document}

\title{Variational Monte Carlo Approach to Partial Differential Equations with Neural Networks}


\author{Moritz Reh}
\email{moritz.reh@kip.uni-heidelberg.de}
\affiliation{Kirchhoff-Institut f\"{u}r Physik, Universit\"{a}t Heidelberg, Im Neuenheimer Feld 227, 69120 Heidelberg, Germany}
\author{Martin G\"{a}rttner}
\email{martin.gaerttner@kip.uni-heidelberg.de}
\affiliation{Kirchhoff-Institut f\"{u}r Physik, Universit\"{a}t Heidelberg, Im Neuenheimer Feld 227, 69120 Heidelberg, Germany}
\affiliation{Physikalisches Institut, Universit\"at Heidelberg, Im Neuenheimer Feld 226, 69120 Heidelberg, Germany}
\affiliation{Institut f\"ur Theoretische Physik, Universit\"at Heidelberg, Philosophenweg 16, 69120 Heidelberg, Germany}

\begin{abstract} 
The accurate numerical solution of partial differential equations (PDEs) is a central task in numerical analysis allowing to model a wide range of natural phenomena by employing specialized solvers depending on the scenario of application.
Here, we develop a variational approach for solving PDEs governing the evolution of high dimensional probability distributions.
Our approach naturally works on the unbounded continuous domain and encodes the full probability density function through its variational parameters, which are adapted dynamically during the evolution to optimally reflect the dynamics of the density. In contrast to previous works, this dynamical adaptation of the parameters is carried out using an explicit prescription avoiding iterative gradient descent.
For the considered benchmark cases we observe excellent agreement with numerical solutions as well as analytical solutions for tasks that are challenging for traditional computational approaches.
\end{abstract}

\maketitle  
\iftrue
\textit{Introduction.} 
The description of nearly all processes in nature is formalized and modelled by means of differential equations, which dictate the evolution of a system given its initial state. Examples include the Navier-Stokes equation in fluid mechanics \cite{Ferziger2002, Cameron2007, Wesseling2001, Spurk2020}, the Schrödinger equation in quantum mechanics \cite{Sakurai2017, Griffiths2018, Schwabl2008}, and the Fokker-Planck equation governing diffusive processes \cite{Kampen2007, Coffey2011, Sornette2001, Freedman1983, Shen2002, Rouse1953, Prakash1999, Reisinger2004}. Analytical solutions of these equations are only available in special cases and, generally, one is forced to resort to numerical techniques. A significant effort during the last century was made to improve the numerical solutions of differential equations \cite{Thomas1995, Hairer2006, Kress1998}. There are numerous properties a numerical solver should ideally fulfill, rendering the field quite diverse, with many specialized solvers being developed \cite{Quarteroni1994}.

Here, we focus on modelling the dynamics of $d$-dimensional probability density functions (PDFs) by means of an ansatz function, which in our case is given by an artificial neural network (ANN), as illustrated in Fig. \ref{fig:IntroFig}. We consider evolution equations of Fokker-Planck form
\begin{equation}
    \label{eq:diff_eq_fundamental}
    \partial_t p = -\sum_i^d\partial_{x_i} \mu_i p + \sum_{ij}^d\partial_{x_i}\partial_{x_j} D_{ij} p,
\end{equation}
where 
$\boldsymbol{\mu}\in\mathbb{R}^d$ is the drift and $D\in\mathbb{R}^{d\times d}$ is the positive semi-definite diffusion matrix and it is understood that $p$, $\boldsymbol{\mu}$ and $D$ are evaluated at position $\mathbf{x}$ and time $t$.

PDFs arise naturally across many disciplines, describing, for example, the phase space evolution of (quantum) matter \cite{Schwabl2006, Zachos2005}, the positions of particles subject to Brownian motion \cite{Freedman1983}, the density of fluids \cite{Ferziger2002} or stock prices in finance \cite{Reisinger2004}. For many of these scenarios the PDF evolution is described by a diffusion process, meaning that the path of a single sampled point evolves according to a stochastic differential equation (SDE) \cite{Tom2015}. In the limit of averaging infinitely many stochastic trajectories one recovers the evolution of the PDF. 

\begin{figure}
    \centering
    \includegraphics[width=\linewidth]{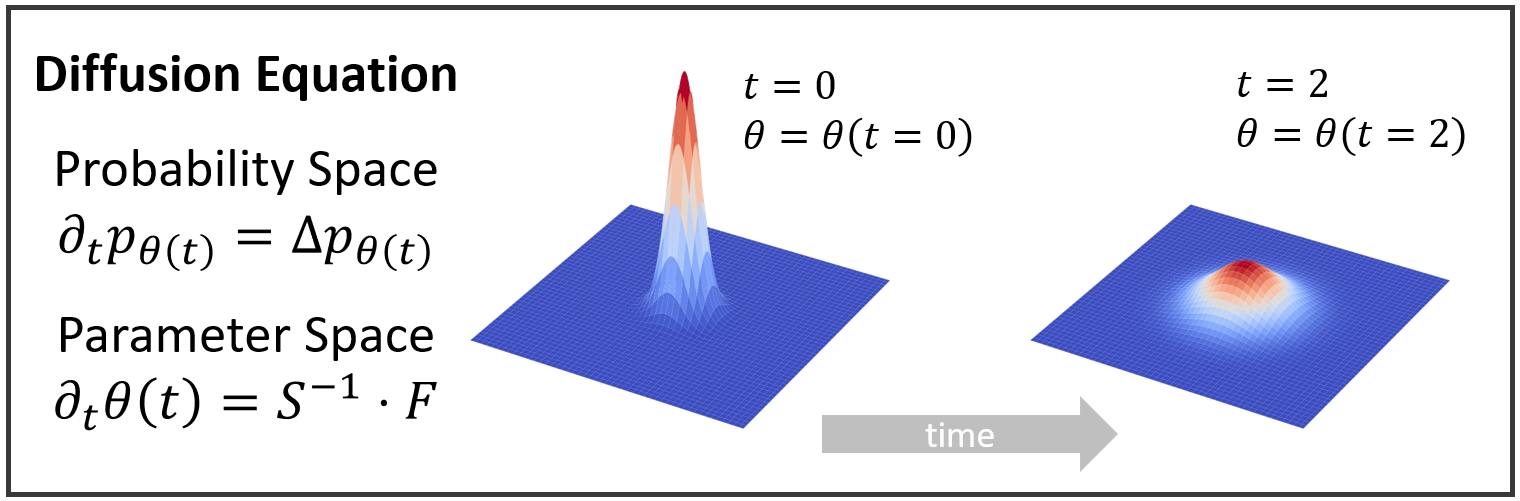}
    \caption{Illustration of the variational approach for a simple diffusion process in 2D. The parameters $\theta(t=0)$ of the artificial neural network encode a gaussian at time $t=0$ and are adapted such that they accurately track the time evolution dictated by the diffusion equation until later times $t=2$, representing a gaussian with increased variance.}
    \label{fig:IntroFig}
\end{figure}

Consequently, the temporal evolution of probability densities can be obtained by either directly solving Eq.~\eqref{eq:diff_eq_fundamental} via spatial discretization (grid based solvers), or by solving the corresponding stochastic dynamics for a large number of sample points (particle based solvers). The former approach, while allowing to control the discretization error via the grid spacing, suffers from the curse of dimensionality \cite{Kang2017, Griebel2005} as the computational cost scales exponentially in the spatial dimension, restricting its applicability to low dimensional cases.
The latter approach solves the SDE associated to the Fokker-Planck equation through the Feynman-Kac formula for an ensemble of points sampled from the initial distribution \cite{DelMoral2004, Oksendal1998}. 
While suited to compute observables, such as moments of the distribution, in high dimensions, there is no direct way to obtain estimates for functionals of the distribution as an expression for $p$ is lacking \cite{Kraskov2004, Singh2016, Ao2022}.

In this work, 
we present a new tool that overcomes the aforementioned limitations of traditional methods by combining
variational Monte-Carlo with normalizing flows (NF). While variational Monte-Carlo is a long established technique in quantum many-body physics \cite{McMillan1965, Ceperley1977, Carleo2017b, Carleo2017}, NFs are a relatively novel class of artificial neural networks also known as invertible neural networks (INNs) \cite{Dinh2016}. They have been applied with remarkable success to long standing problems in statistical physics \cite{Noe2019}, inference and data generation \cite{Papamakarios2019, Dinh2016, Grathwohl2018, Dinh2014, Kingma2018, Ardizzone2018}, as well as quantum field theories \cite{Albergo2021, Pawlowski2022}. Here, we understand the NF as an ansatz function for the time-dependent density.
The choice of the ansatz-function is a degree of freedom in our approach and can be adapted to the problem at hand exploiting prior knowledge about the function class the time-dependent density belongs to.
Among the possible choices, artificial neural networks are a promising class of ansatz functions, as
they may become \textit{universal} function approximators in the infinite parameter limit, which applies to lesser extent to NFs \cite{Kong2020, Teshima2020}.
Adjusting the parameters of the ansatz function to the dynamics dictated by Eq.~\eqref{eq:diff_eq_fundamental} is achieved by a time-dependent variational principle (TDVP), which maps the dynamics of the PDF onto the variational manifold generated by the ansatz function \cite{Carleo2017b, Carleo2017, Reh2021}. 
Crucially, the approach is self-contained and at no point relies on data generated from other solvers, in contrast to prior works using neural networks to solve PDEs \cite{Brandstetter2022, Lu2021, Li2020, Beck2021}, allowing us to obtain numerical solutions for tasks that are challenging for grid-based or particle-based solvers. 
Our approach differs from the popular physics informed neural networks (PINN) \cite{Raissi2019} 
and \cite{Feng2022}
in that we do not carry out a costly global gradient-descent based optimization in each time step to update the models' parameters, but rather follow an explicit, analytically derived time derivative of the network parameters which is given by the TDVP.
We are particularly interested in high-dimensional scenarios
which are infeasible to solve with grid-based methods and in quantities which are not easily obtainable by modelling many stochastic processes, such as
functionals of the PDF.
Indeed we show that, 
using the developed approach,
we can reliably estimate differential entropies
in a Monte Carlo fashion requiring only a few thousand samples.
We benchmark our approach for the case of an eight-dimensional heat equation and a six-dimensional dissipative phase space evolution.

\textit{Normalizing flows.}
While we employ neural networks as ansatz functions, we emphasize that the derived TDVP is applicable to any parameterized density, such as Gaussian mixture models or energy-based estimators.
We use NFs \cite{Dinh2016, Papamakarios2019} to model densities as they
have many desirable properties, among which are (i) a guarantee of normalization for any set of parameters $\theta$, (ii) a tractable likelihood and (iii) the ability to generate independent samples without the need to resort to Markov Chains. NFs parameterize densities by assuming a latent distribution $\pi$ which is transformed into the distribution of interest by a trainable and invertible map $\mathbf{f}_\theta$,
\begin{equation}
    \mathbf{x} = \mathbf{f}_\theta(\mathbf{z}) \mbox{ with } \mathbf{z}\sim\pi.
\end{equation}
Usually, $\pi$ is chosen to be a `simple' distribution, e.g. a Gaussian, such that its samples $\mathbf{z}$ can be generated easily.
The probability associated with the point $\mathbf{x}$ is proportional to $\pi(\mathbf{f}_\theta^{-1}(\mathbf{x}))$ times the determinant of the Jacobian of the transformation,
\begin{equation}
    p_\theta(\mathbf{x}) = \pi(\textbf{f}_\theta^{-1}(\mathbf{x})) \biggm\lvert \det\left(\frac{\partial \mathbf{f}_\theta^{-1}(\mathbf{x})}{\partial \mathbf{x}}\right) \biggm\lvert.
\end{equation}

The function $\mathbf{f}_\theta$ is composed from a series of invertible transformations $\mathbf{f}_\theta=\boldsymbol{\varphi}^1_\theta\circ..\circ\boldsymbol{\varphi}^N_\theta$ which are explained in detail in the Supplemental Material (SM) \cite{SM}. Importantly, the Jacobian of the function is tractable meaning that its determinant is efficiently inferred when computing a forward pass, an operation carried out whenever the real space probability is evaluated at some point of interest.
By stacking many of these `coupling blocks' $\boldsymbol{\varphi}^i$, the function $\mathbf{f}_\theta$ becomes an expressive coordinate transform, that is, however, incapable of changing the tail behavior of the latent space distribution \cite{Jaini2019}.
We overcome this problem by dynamically adapting the latent space distribution $\pi$ to reflect dynamical changes in the tails of the distribution. This is explained in more detail below and in the SM \cite{SM}.

\textit{Time-Dependent Variational Principle.} The idea of the TDVP originated in the context of Variational Monte Carlo (VMC) \cite{McMillan1965} where it has been applied extensively to solve problems in quantum-many-body physics, 
with a growing interest in the use of neural networks as variational ansatz functions \cite{Carleo2017b, Carleo2017, Schmitt2020, Reh2021}.
Its aim is to locally search for the closest approximation to the dynamics of the density within the variational manifold. Concretely, one aims to solve
\begin{equation}
    \label{eq:tdvp_cond}
    \underset{\dot\theta}{\mathrm{argmin}} \: \mathcal{D}(p_{\theta(t) + \tau \dot\theta}, p_{\theta(t)} + \tau \dot{p}_{\theta(t)})
\end{equation}
where $\mathcal{D}$ is a suitable distance measure between probability distributions, $\tau$ denotes a small time step, $\dot{p}_{\theta(t)}$ is the derivative given by Eq.~\eqref{eq:diff_eq_fundamental}, and $\dot{\theta}$ is the unknown corresponding parameter time derivative. 
The solution to Eq.~\eqref{eq:tdvp_cond} can be found by requiring the derivative with respect to $\dot{\theta}$ to be zero. 
By expanding Eq.~\eqref{eq:tdvp_cond} to second order in $\tau$ one finds
\begin{equation}
    \label{eq:tdvp_Solution}
    S_{kk'} \dot{\theta}_{k'} = F_k.
\end{equation}
We defer the details of this derivation to the SM \cite{SM}.
Here $S_{kk'}=\langle O_k(\textbf{x}) O_{k'}(\textbf{x})\rangle_{\mathbf{x}\sim p_\theta(t)}$ denotes the Fisher information metric and $F_k = \langle O_k(\textbf{x}) \partial_t\log (p_{\theta(t)}(\textbf{x}))\rangle_{\mathbf{x}\sim p_{\theta(t)}}$ is a force term, where $O_k$ denotes the (logarithmic) variational derivative $O_k(\textbf{x}) = \partial_{\theta_k} \log (p_{\theta(t)}(\textbf{x}))$ and $\partial_t\log (p_{\theta(t)}(\textbf{x}))$ is given by the RHS of the to be solved PDE. Here $\langle \cdot\rangle_{\mathbf{x}\sim p_\theta(t)}$ denotes an expectation value evaluated through Monte Carlo sampling from the model distribution $p_\theta(t)$.
Notice, that we heavily rely on the differentiability of the ansatz function $p_{\theta(t)}$ with respect to both variational parameters and spatial coordinates.
The latter frequently appear on the RHS of Eq.~\eqref{eq:diff_eq_fundamental} and are thus required for computing $\dot{p}_{\theta(t)}$.
This is in striking contrast to grid-based techniques which require making grid cells finer for higher accuracy. Here, instead, we have access to the exact derivatives through automatic differentiation. The choice of distance measure to compare the two probability distributions is not arbitrary as the form of $S$ and $F$ directly depends on it. In order to obtain expressions of $S$ and $F$ that can be efficiently estimated through a finite number of samples, we found that both the Hellinger distance $\mathcal{D}_H(p, q) = 1 - F(p, q) = 1 - \int \sqrt{pq}d\mathbf{x}$ and the Kullback-Leibler (KL) divergence $\mathcal{D}_{KL}(p, q) = \int p\log(p/q) d\mathbf{x}$ yield the same result of the desired form. 
Care has to be taken when solving Eq.~\eqref{eq:tdvp_Solution} for $\dot\theta$, as the inverse of $S$ may not exist.
This is the case if directions in parameter space are present along which the probabilities are stationary, which can be dealt with by regularization procedures \cite{Carleo2017, Schmitt2020}.

\textit{Problem Setup.}
We are interested in solving initial value problems, for which the initial density distribution $p(0, \mathbf{x}) = u(\mathbf{x})$ is given along with the RHS of Eq.~\eqref{eq:diff_eq_fundamental} which governs its evolution.
To exactly encode the initial distribution $u(\mathbf{x})$ in the model $p_{\theta(t=0)}$, the latent distribution is set to $u(\mathbf{x})$ and the parameters of the map $\mathbf{f}_{\theta(t=0)}$ are chosen such that it represents the identity map $\mathbf{f}_{\theta(t=0)}(\mathbf{x}) = \mathbf{x}$. If the initial distribution cannot be given in closed form and therefore cannot be set analytically as the latent space distribution $\pi$, the network may be trained on its samples to approximately encode it at time $t=0$. Then a solver is used which integrates the parameters according to Eq.~\eqref{eq:tdvp_Solution}.

\textit{Application 1: Diffusion in High Dimensions.}
As a first benchmark scenario we consider the heat equation in $d=8$ dimensions. The heat equation appears across many disciplines ranging from engineering \cite{Poirier2016, Annaratone2010} and molecular motion \cite{Freedman1983} to the pricing of financial derivatives given by the famous Black-Scholes equation \cite{Black1973, Guillaume2018} and
reads
\begin{equation}
    \label{eq:heat_eq}
    \partial_t p(t, \mathbf{x}) = D\Delta_\mathbf{x} p(t, \mathbf{x}).
\end{equation}
Importantly, an analytical solution exists against which we can benchmark, making the described scenario a good showcase of the proposed approach. The solution is given by a convolution of the initial distribution $p(0, \mathbf{x})$ with the `heat kernel' $\Phi(t, \mathbf{x})=(4\pi t)^{-(d/2)}\exp(-\mathbf{x}^2/4 Dt)$ \cite{Davies1989}, which is the Green's function to Eq.~\eqref{eq:heat_eq}, such that
\begin{equation}
    p(t, \mathbf{x}) = \int p(0, \mathbf{y}) \Phi(t, \mathbf{x}-\mathbf{y})d\mathbf{y}.
\end{equation}
We aim to observe the growth of the differential entropy 
\begin{equation}
    \begin{split}
    \label{eqn:diffEntropy}
    S(t) &= -\int p(t, \mathbf{x})\log(p(t, \mathbf{x}))d\mathbf{x}\\
    &=-\langle\log(p(t, \mathbf{x}))\rangle_{\mathbf{x}\sim p(t, \mathbf{x})}    
    \end{split}
\end{equation}
with time, a task, which is challenging or even intractable using other numerical techniques in high dimensions for the reasons mentioned above \cite{Kang2017, Griebel2005, Kraskov2004, Singh2016, Ao2022}. 
In the case of a Gaussian distribution for $p(0, \mathbf{x})$ with zero mean and unit covariance matrix, we obtain a Gaussian of larger variance at later points in time, in which case we observe perfect agreement between the analytical solution and the one obtained using the INN as shown in Fig.~\ref{fig:secondFig}. 
If we choose a Student-$t$ distribution as initial distribution, i.e.
\begin{equation}
    p(0, \mathbf{x})\propto \left(1 + \frac{\mathbf{x}^2}{\nu}\right)^{-(\nu + d) / 2}
\end{equation}
with $\nu=2$ we can no longer compare to the analytical solution as the involved integrals become infeasible to solve. However, by exploiting the spherical symmetry of the problem, we can map the evolution to an effective 1D problem of the radial dependency of $p$ which we can approximately solve on a grid using finite differences.
The grid based solution and that obtained using the INN are generally in good agreement. We observe a slight difference which we attribute to technical challenges of the grid-based approach, which we discuss more elaborately in the SM \cite{SM}.

\begin{figure}
    \centering
    \includegraphics[width=\linewidth]{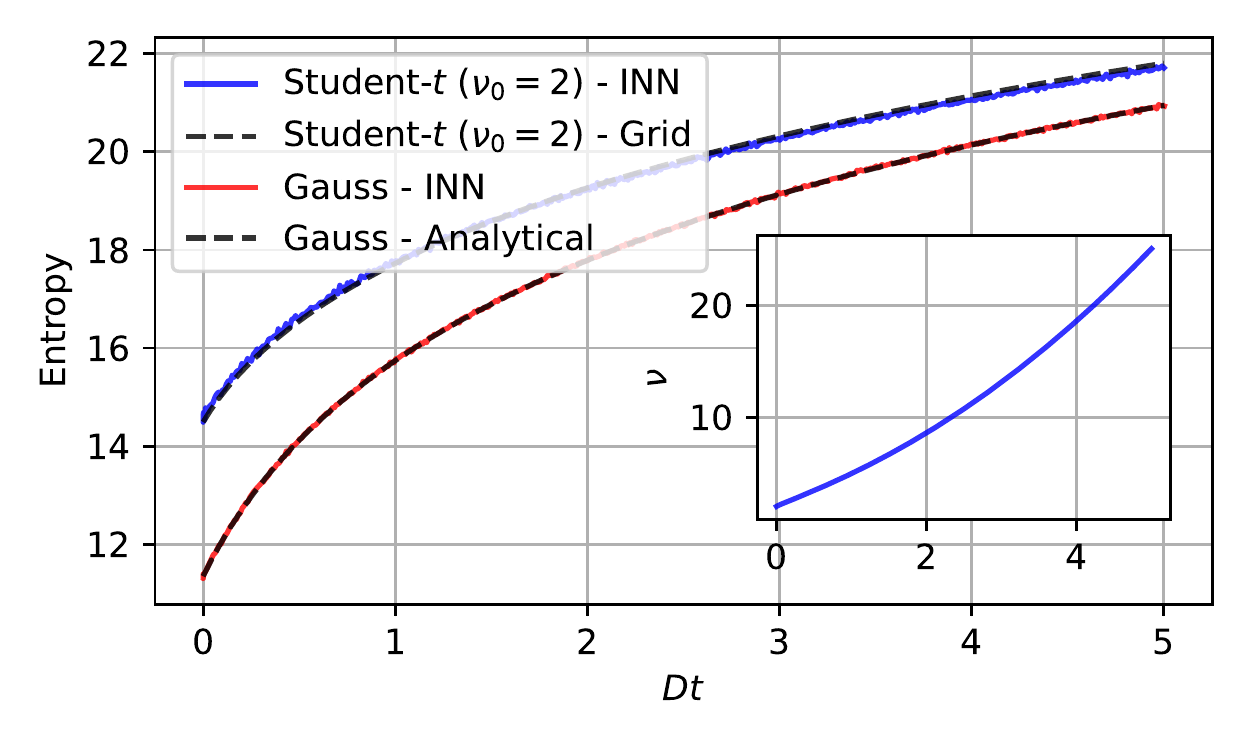}
    \caption{Evolution of the differential entropy under a heat equation for different initial distributions. 
    While analytical comparison data is available in the case of a Gaussian initial distribution, we compare to numerical data for the Student-$t$ obtained with finite differences on a 1D Grid.
    Inset: Adjustment of the latent space distribution $\pi$ by changing its parameter $\nu$ dynamically in time.
    }
    \label{fig:secondFig}
\end{figure}

\begin{figure*}
    \centering
    \includegraphics[width=\textwidth]{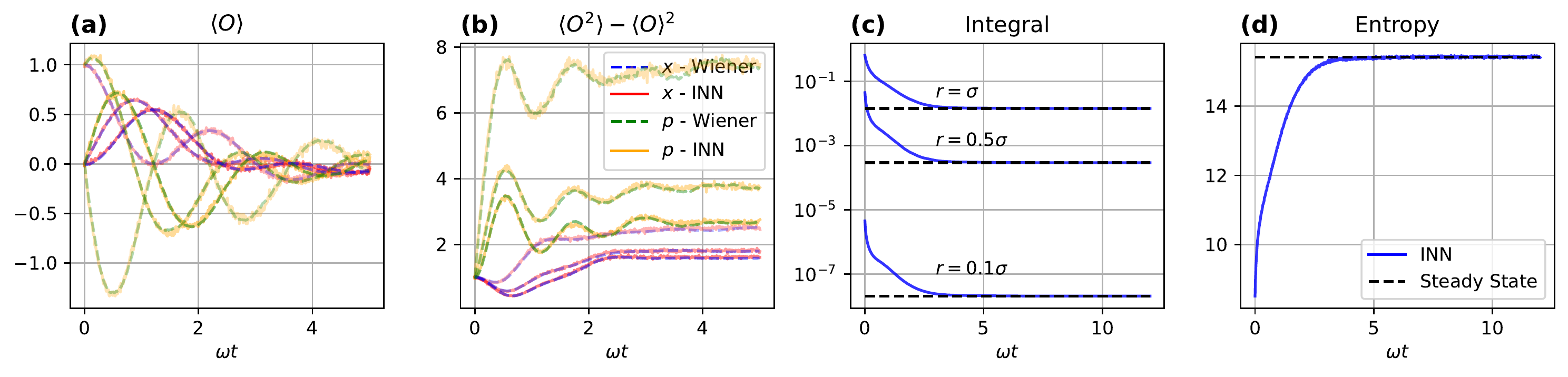}
    \caption{(a) and (b): Evolution of the first two moments of the phase space distribution $\rho$ estimated from 10.000 samples for three coupled harmonic oscillators with dissipation given by the temperatures $k_BT /m\omega^2 =(10, 3, 1)$ and all other parameters chosen to be unity. The initial distribution in phase space is a Gaussian with unit variance centered at the position $\mathbf{x}=(1, 0, 0)^T$ and momentum $\mathbf{p}=(0, 1, 0)^T$. 
    (c) and (d): Three uncoupled oscillators ($k=0$) coupled to the same heat bath at temperature $k_BT/m\omega^2=10$. 
    In (c) the value of the six-dimensional integral around a hypersphere with radius $r$ centered at the origin is shown, while in (d) the estimation of the differential entropy Eq.~\eqref{eqn:diffEntropy} using the INN is displayed. Both are shown to converge to the expected value of the steady state. The initial distribution is a Gaussian with identity covariance matrix centered at $\mathbf{x}=(1, 1, 1)^T$ and zero momentum.}
    \label{fig:thirdFig}
\end{figure*}

\textit{Application 2: Diffusion in classical phase space.} 
As a second demonstration of the proposed approach we consider classical Hamiltonian dynamics in phase space with additional diffusion. 
Concretely, we choose the Hamiltonian $H$ to represent coupled harmonic oscillators (coupling strength $k$) which are in contact with heat baths of different temperatures $T_i$, such that the solution does not factorize in the eigenbasis of $H$. We provide the Hamiltonian and its generated phase space flow in the SM \cite{SM}.
The heat baths lead to a diffusion in phase space, which implies that sampled points of the distribution evolve according to an SDE.
We show that the INN faithfully estimates moments of the distribution, probabilities (i.e. integrals over finite domains) as well as functionals of the PDF that correspond to integrals over the entire domain. 


The described system obeys the following Fokker-Planck equation \cite{Presilla1997}
\begin{equation}
\begin{split}
\partial_t \rho(t, \mathbf{x}, \mathbf{p}) = &\left[-\partial_\mathbf{p}H \cdot \partial_\mathbf{x} + \partial_\mathbf{x}H\cdot\partial_\mathbf{p}+ \right.\\
&\left. \gamma\left(\mathbf{p}\cdot\partial_\mathbf{p} + m k_B \mathsmaller{\sum}_i T_i\partial_{p_i}^2 \right)\right] \rho(t, \mathbf{x}, \mathbf{p}),    
\end{split}
\end{equation}
whose corresponding stochastic differential equation is given by \cite{Presilla1997}
\begin{equation}
\begin{split}
    & dx_i = \partial_{p_i}H dt,\\
    & dp_i = -\left[\gamma p_i + \partial_{x_i}H\right]dt + \sqrt{\lambda_i} dw_i.
\end{split}
\end{equation}
Here, $\lambda_i=\sqrt{2m\gamma k_B T_i}$, $dw_i=\Pi_i \sqrt{dt}$ is the Wiener process with zero average $\langle dw_i\rangle=0$ and standard scaling $\langle dw_i^2 \rangle=dt$ implying that $\Pi$ is drawn from a standard Gaussian $\Pi \sim \mathcal{N}(0, 1)$.
For simplicity we choose all quantities except $T_i$ equal to unity. 

In the case of heat baths of equal temperatures $T_i=T$ and vanishing coupling ($k=0$) the system, in contrast, assumes a thermal steady state of Gaussian form in the long time limit given by the Gibbs-Ensemble
\begin{equation}
\begin{split}
    \label{eqn:gauss_steady_state}
    \rho_{SS} & = \exp(-H/k_BT) / Z \\
             & = \exp(-\frac{1}{2}(m\omega^2\mathbf{x}^2 + \mathbf{p}^2/m)/k_BT) / Z,
\end{split}
\end{equation}
with $Z=\int \exp(-H/k_BT) d\mathbf{x}d\mathbf{p}$ the partition function, where the Gaussian form allows to compare against analytical results.

We consider four quantities of interest which we evaluate by drawing 10.000 samples from the INN, see Fig.~\ref{fig:thirdFig}. The first two quantities are the means and variances of the distribution evolved for the case of different $T_i$ and $k=1$. Here, comparison against estimates from solving the SDE for the same number of sampled points is straight forward and one observes excellent agreement between both methods.
To obtain an easy benchmark case for integral and entropy estimation, we choose $k=0$ and $T_i=T$ such that the steady state is Gaussian, see Eq.~\eqref{eqn:gauss_steady_state}.
We choose the integration volumes to be hyperspheres of radius $r$ centered at the origin allowing for analytical evaluation of the Gaussian integral.
The values of these integrals correspond to the probability of finding the system inside the hypersphere.
Using the INN, we can estimate such integrals in a Monte-Carlo fashion by uniformly sampling points $\mathbf{x}_i$ from inside the integration domain and average the associated probabilities $p_\theta(\mathbf{x}_i)$, which are shown to converge to the analytically obtained steady-state value in Fig.~\ref{fig:secondFig}(c).

Finally, we again focus on the differential entropy (Eq.~\eqref{eqn:diffEntropy}), where Fig.~\ref{fig:thirdFig}(d) shows that our method succeeds to predict the differential entropy with low noise while converging to the expected steady state value.

\textit{Conclusion and Outlook.}
We have introduced a variational approach to the dynamics of continuous probability distributions using normalizing flows and demonstrated its power by applying it to paradigmatic benchmark problems. Our method is widely applicable, even beyond the Fokker-Planck form \eqref{eq:diff_eq_fundamental}, e.g. to cases with non-local terms \cite{Lin2021a}.
Its unique strength lies in estimating functionals of probability densities in high dimensions enabled by the availability of exact samples with tractable likelihood.
We emphasize that other approaches such as PINN \cite{Raissi2019} require solving a large-scale non-convex optimization problem in each time step, which the TDVP replaces by the explicit update rule \eqref{eq:tdvp_Solution} (see \cite{SM} for further discussion).
The form of the ansatz function can be chosen flexibly and is not required to be a neural network. The only restrictions are that (i) samples from its distribution may be obtained and (ii) derivatives with respect to inputs and parameters are computable. While building normalizing flows using stacked coupling blocks is a popular approach, other flow architectures exist and it would be interesting to investigate their potential in solving PDEs in the future. Since the TDVP can also work with non-normalized probabilities, also energy based models would be viable ansatz functions although this would mean that samples would have to be obtained by resorting to Markov-Chains.

For the utilized architecture we found that challenges exist when trying to solve chaotic dynamics. We believe this to be caused by the high amount of information of the phase space distribution which needs to be encoded using comparably few parameters. 
Additionally, we found it challenging to model distributions whose tail behaviour deviated from that of the latent space distribution. In the example shown in Fig.~\ref{fig:secondFig} this could be dealt with by elevating $\nu$ to be a variational parameter, which would tend to infinity for late times, representing the exact tail behaviour of the real space distribution. 
However, if the real space tail behaviour cannot be accurately modelled in latent space,
e.g. because its form is not known beforehand, one cannot expect to accurately model the distribution on the entire domain.


\textit{Code \& Data availability.}
The code used for this project is based on the \texttt{jVMC} library \cite{jvmc}, making use of \texttt{flax} \cite{flax2020github} and \texttt{jax} \cite{Bradbury2018} and is available under \href{https://github.com/RehMoritz/vmc_pde}{GitHub: RehMoritz/vmc\_pde}. The repository also contains the data from Fig.~\ref{fig:secondFig} and Fig.~\ref{fig:thirdFig}.

\acknowledgments
\textit{Acknowledgments.}
The authors thank Markus Schmitt, Julian Urban, Ullrich Köthe and Peter Sorrenson for helpful discussions. 
This work is supported by the Deutsche Forschungsgemeinschaft (DFG, German Research Foundation) under Germany’s Excellence Strategy EXC2181/1-390900948 (the Heidelberg STRUCTURES Excellence Cluster) and within the Collaborative Research Center SFB1225 (ISOQUANT). This work was partially financed by the Baden-Württemberg Stiftung gGmbH. The authors acknowledge support by the state of Baden-Württemberg through bwHPC
and the German Research Foundation (DFG) through Grant No INST 40/575-1 FUGG (JUSTUS 2 cluster). The authors gratefully acknowledge the Gauss Centre for Supercomputing e.V. (www.gauss-centre.eu) for funding this project by providing computing time through the John von Neumann Institute for Computing (NIC) on the GCS Supercomputer JUWELS \cite{JUWELS} at Jülich Supercomputing Centre (JSC).

\bibliography{refs}

\begin{thebibliography}{65}%
\makeatletter
\providecommand \@ifxundefined [1]{%
 \@ifx{#1\undefined}
}%
\providecommand \@ifnum [1]{%
 \ifnum #1\expandafter \@firstoftwo
 \else \expandafter \@secondoftwo
 \fi
}%
\providecommand \@ifx [1]{%
 \ifx #1\expandafter \@firstoftwo
 \else \expandafter \@secondoftwo
 \fi
}%
\providecommand \natexlab [1]{#1}%
\providecommand \enquote  [1]{``#1''}%
\providecommand \bibnamefont  [1]{#1}%
\providecommand \bibfnamefont [1]{#1}%
\providecommand \citenamefont [1]{#1}%
\providecommand \href@noop [0]{\@secondoftwo}%
\providecommand \href [0]{\begingroup \@sanitize@url \@href}%
\providecommand \@href[1]{\@@startlink{#1}\@@href}%
\providecommand \@@href[1]{\endgroup#1\@@endlink}%
\providecommand \@sanitize@url [0]{\catcode `\\12\catcode `\$12\catcode
  `\&12\catcode `\#12\catcode `\^12\catcode `\_12\catcode `\%12\relax}%
\providecommand \@@startlink[1]{}%
\providecommand \@@endlink[0]{}%
\providecommand \url  [0]{\begingroup\@sanitize@url \@url }%
\providecommand \@url [1]{\endgroup\@href {#1}{\urlprefix }}%
\providecommand \urlprefix  [0]{URL }%
\providecommand \Eprint [0]{\href }%
\providecommand \doibase [0]{http://dx.doi.org/}%
\providecommand \selectlanguage [0]{\@gobble}%
\providecommand \bibinfo  [0]{\@secondoftwo}%
\providecommand \bibfield  [0]{\@secondoftwo}%
\providecommand \translation [1]{[#1]}%
\providecommand \BibitemOpen [0]{}%
\providecommand \bibitemStop [0]{}%
\providecommand \bibitemNoStop [0]{.\EOS\space}%
\providecommand \EOS [0]{\spacefactor3000\relax}%
\providecommand \BibitemShut  [1]{\csname bibitem#1\endcsname}%
\let\auto@bib@innerbib\@empty
\bibitem [{\citenamefont {Ferziger}\ and\ \citenamefont
  {Peri{\'{c}}}(2002)}]{Ferziger2002}%
  \BibitemOpen
  \bibfield  {author} {\bibinfo {author} {\bibfnamefont {J.~H.}\ \bibnamefont
  {Ferziger}}\ and\ \bibinfo {author} {\bibfnamefont {M.}~\bibnamefont
  {Peri{\'{c}}}},\ }\href {\doibase https://doi.org/10.1007/978-3-642-56026-2}
  {\emph {\bibinfo {title} {Computational Methods for Fluid Dynamics}}}\
  (\bibinfo  {publisher} {Springer Berlin Heidelberg},\ \bibinfo {year}
  {2002})\BibitemShut {NoStop}%
\bibitem [{\citenamefont {Tropea}\ \emph {et~al.}(2007)\citenamefont {Tropea},
  \citenamefont {Yarin},\ and\ \citenamefont {Foss}}]{Cameron2007}%
  \BibitemOpen
  \bibfield  {author} {\bibinfo {author} {\bibfnamefont {C.}~\bibnamefont
  {Tropea}}, \bibinfo {author} {\bibfnamefont {A.}~\bibnamefont {Yarin}}, \
  and\ \bibinfo {author} {\bibfnamefont {J.}~\bibnamefont {Foss}},\ }\href
  {\doibase 10.1007/978-3-540-30299-5} {\emph {\bibinfo {title} {Springer
  Handbook of Experimental Fluid Mechanics}}}\ (\bibinfo  {publisher} {Springer
  Berlin Heidelberg},\ \bibinfo {year} {2007})\BibitemShut {NoStop}%
\bibitem [{\citenamefont {Wesseling}(2001)}]{Wesseling2001}%
  \BibitemOpen
  \bibfield  {author} {\bibinfo {author} {\bibfnamefont {P.}~\bibnamefont
  {Wesseling}},\ }\href {\doibase https://doi.org/10.1007/978-3-642-05146-3}
  {\emph {\bibinfo {title} {Principles of Computational Fluid Dynamics}}}\
  (\bibinfo  {publisher} {Springer Berlin Heidelberg},\ \bibinfo {year}
  {2001})\BibitemShut {NoStop}%
\bibitem [{\citenamefont {Spurk}\ and\ \citenamefont
  {Aksel}(2020)}]{Spurk2020}%
  \BibitemOpen
  \bibfield  {author} {\bibinfo {author} {\bibfnamefont {J.~H.}\ \bibnamefont
  {Spurk}}\ and\ \bibinfo {author} {\bibfnamefont {N.}~\bibnamefont {Aksel}},\
  }\href {\doibase https://doi.org/10.1007/978-3-030-30259-7} {\emph {\bibinfo
  {title} {Fluid Mechanics}}}\ (\bibinfo  {publisher} {Springer International
  Publishing},\ \bibinfo {year} {2020})\BibitemShut {NoStop}%
\bibitem [{\citenamefont {Sakurai}\ and\ \citenamefont
  {Napolitano}(2017)}]{Sakurai2017}%
  \BibitemOpen
  \bibfield  {author} {\bibinfo {author} {\bibfnamefont {J.~J.}\ \bibnamefont
  {Sakurai}}\ and\ \bibinfo {author} {\bibfnamefont {J.}~\bibnamefont
  {Napolitano}},\ }\href {\doibase https://doi.org/10.1017/9781108499996}
  {\emph {\bibinfo {title} {Modern Quantum Mechanics}}}\ (\bibinfo  {publisher}
  {Cambridge University Press},\ \bibinfo {year} {2017})\BibitemShut {NoStop}%
\bibitem [{\citenamefont {Griffiths}\ and\ \citenamefont
  {Schroeter}(2018)}]{Griffiths2018}%
  \BibitemOpen
  \bibfield  {author} {\bibinfo {author} {\bibfnamefont {D.~J.}\ \bibnamefont
  {Griffiths}}\ and\ \bibinfo {author} {\bibfnamefont {D.~F.}\ \bibnamefont
  {Schroeter}},\ }\href {\doibase https://doi.org/10.1017/9781316995433} {\emph
  {\bibinfo {title} {Introduction to Quantum Mechanics}}}\ (\bibinfo
  {publisher} {Cambridge University Press},\ \bibinfo {year}
  {2018})\BibitemShut {NoStop}%
\bibitem [{\citenamefont {Schwabl}(2008)}]{Schwabl2008}%
  \BibitemOpen
  \bibfield  {author} {\bibinfo {author} {\bibfnamefont {F.}~\bibnamefont
  {Schwabl}},\ }\href {\doibase 10.1007/978-3-540-85062-5} {\emph {\bibinfo
  {title} {Advanced Quantum Mechanics}}}\ (\bibinfo  {publisher} {Springer
  Berlin Heidelberg},\ \bibinfo {year} {2008})\BibitemShut {NoStop}%
\bibitem [{\citenamefont {Kampen}(2007)}]{Kampen2007}%
  \BibitemOpen
  \bibfield  {author} {\bibinfo {author} {\bibfnamefont {N.~V.}\ \bibnamefont
  {Kampen}},\ }\href {\doibase
  https://doi.org/10.1016/B978-0-444-52965-7.X5000-4} {\emph {\bibinfo {title}
  {Stochastic Processes in Physics and Chemistry}}}\ (\bibinfo  {publisher}
  {Elsevier},\ \bibinfo {year} {2007})\BibitemShut {NoStop}%
\bibitem [{\citenamefont {Coffey}\ and\ \citenamefont
  {Kalmykov}(2011)}]{Coffey2011}%
  \BibitemOpen
  \bibfield  {author} {\bibinfo {author} {\bibfnamefont {W.~T.}\ \bibnamefont
  {Coffey}}\ and\ \bibinfo {author} {\bibfnamefont {Y.~P.}\ \bibnamefont
  {Kalmykov}},\ }\href {\doibase https://doi.org/10.1142/8195} {\emph {\bibinfo
  {title} {The Langevin Equation}}}\ (\bibinfo  {publisher} {World Scientific
  Publishing Company},\ \bibinfo {year} {2011})\BibitemShut {NoStop}%
\bibitem [{\citenamefont {Sornette}(2001)}]{Sornette2001}%
  \BibitemOpen
  \bibfield  {author} {\bibinfo {author} {\bibfnamefont {D.}~\bibnamefont
  {Sornette}},\ }\href {\doibase https://doi.org/10.1016/S0378-4371(00)00571-9}
  {\bibfield  {journal} {\bibinfo  {journal} {Physica A: Statistical Mechanics
  and its Applications}\ }\textbf {\bibinfo {volume} {290}},\ \bibinfo {pages}
  {211} (\bibinfo {year} {2001})}\BibitemShut {NoStop}%
\bibitem [{\citenamefont {Freedman}(1983)}]{Freedman1983}%
  \BibitemOpen
  \bibfield  {author} {\bibinfo {author} {\bibfnamefont {D.}~\bibnamefont
  {Freedman}},\ }\href {\doibase https://doi.org/10.1007/978-1-4615-6574-1}
  {\emph {\bibinfo {title} {Brownian Motion and Diffusion}}}\ (\bibinfo
  {publisher} {Springer New York},\ \bibinfo {year} {1983})\BibitemShut
  {NoStop}%
\bibitem [{\citenamefont {Shen}\ \emph {et~al.}(2002)\citenamefont {Shen},
  \citenamefont {Chen}, \citenamefont {Dai},\ and\ \citenamefont
  {Dai}}]{Shen2002}%
  \BibitemOpen
  \bibfield  {author} {\bibinfo {author} {\bibfnamefont {X.}~\bibnamefont
  {Shen}}, \bibinfo {author} {\bibfnamefont {H.}~\bibnamefont {Chen}}, \bibinfo
  {author} {\bibfnamefont {J.}~\bibnamefont {Dai}}, \ and\ \bibinfo {author}
  {\bibfnamefont {W.}~\bibnamefont {Dai}},\ }\href {\doibase
  https://doi.org/10.1023/A:1019942711261} {\bibfield  {journal} {\bibinfo
  {journal} {Queueing Systems}\ }\textbf {\bibinfo {volume} {42}},\ \bibinfo
  {pages} {33} (\bibinfo {year} {2002})}\BibitemShut {NoStop}%
\bibitem [{\citenamefont {Rouse}(1953)}]{Rouse1953}%
  \BibitemOpen
  \bibfield  {author} {\bibinfo {author} {\bibfnamefont {P.~E.}\ \bibnamefont
  {Rouse}},\ }\href {\doibase https://doi.org/10.1063/1.1699180} {\bibfield
  {journal} {\bibinfo  {journal} {The Journal of Chemical Physics}\ }\textbf
  {\bibinfo {volume} {21}},\ \bibinfo {pages} {1272} (\bibinfo {year}
  {1953})}\BibitemShut {NoStop}%
\bibitem [{\citenamefont {Prakash}\ and\ \citenamefont
  {Öttinger}(1999)}]{Prakash1999}%
  \BibitemOpen
  \bibfield  {author} {\bibinfo {author} {\bibfnamefont {J.~R.}\ \bibnamefont
  {Prakash}}\ and\ \bibinfo {author} {\bibfnamefont {H.~C.}\ \bibnamefont
  {Öttinger}},\ }\href {\doibase https://doi.org/10.1021/ma981534b} {\bibfield
   {journal} {\bibinfo  {journal} {Macromolecules}\ }\textbf {\bibinfo {volume}
  {32}},\ \bibinfo {pages} {2028} (\bibinfo {year} {1999})}\BibitemShut
  {NoStop}%
\bibitem [{\citenamefont {Reisinger}(2004)}]{Reisinger2004}%
  \BibitemOpen
  \bibfield  {author} {\bibinfo {author} {\bibfnamefont {C.}~\bibnamefont
  {Reisinger}},\ }\href {\doibase 10.11588/heidok.00004954} {\emph {\bibinfo
  {title} {Numerische Methoden für hochdimensionale parabolische Gleichungen
  am Beispiel von Optionspreisaufgaben}}}\ (\bibinfo  {publisher} {Heidelberg
  University Library},\ \bibinfo {year} {2004})\BibitemShut {NoStop}%
\bibitem [{\citenamefont {Thomas}(1995)}]{Thomas1995}%
  \BibitemOpen
  \bibfield  {author} {\bibinfo {author} {\bibfnamefont {J.~W.}\ \bibnamefont
  {Thomas}},\ }\href {\doibase https://doi.org/10.1007/978-1-4899-7278-1}
  {\emph {\bibinfo {title} {Numerical Partial Differential Equations: Finite
  Difference Methods}}}\ (\bibinfo  {publisher} {Springer New York},\ \bibinfo
  {year} {1995})\BibitemShut {NoStop}%
\bibitem [{\citenamefont {Hairer}(2006)}]{Hairer2006}%
  \BibitemOpen
  \bibfield  {author} {\bibinfo {author} {\bibfnamefont {E.}~\bibnamefont
  {Hairer}},\ }\href {\doibase https://doi.org/10.1007/3-540-30666-8} {\emph
  {\bibinfo {title} {Geometric Numerical Integration}}}\ (\bibinfo  {publisher}
  {Springer-Verlag},\ \bibinfo {year} {2006})\BibitemShut {NoStop}%
\bibitem [{\citenamefont {Kress}(1998)}]{Kress1998}%
  \BibitemOpen
  \bibfield  {author} {\bibinfo {author} {\bibfnamefont {R.}~\bibnamefont
  {Kress}},\ }\href {\doibase https://doi.org/10.1007/978-1-4612-0599-9_9}
  {\emph {\bibinfo {title} {Graduate Texts in Mathematics}}}\ (\bibinfo
  {publisher} {Springer New York},\ \bibinfo {year} {1998})\ pp.\ \bibinfo
  {pages} {189--224}\BibitemShut {NoStop}%
\bibitem [{\citenamefont {Quarteroni}\ and\ \citenamefont
  {Valli}(1994)}]{Quarteroni1994}%
  \BibitemOpen
  \bibfield  {author} {\bibinfo {author} {\bibfnamefont {A.}~\bibnamefont
  {Quarteroni}}\ and\ \bibinfo {author} {\bibfnamefont {A.}~\bibnamefont
  {Valli}},\ }\href {\doibase https://doi.org/10.1007/978-3-540-85268-1} {\emph
  {\bibinfo {title} {Numerical Approximation of Partial Differential
  Equations}}}\ (\bibinfo  {publisher} {Springer Berlin Heidelberg},\ \bibinfo
  {year} {1994})\BibitemShut {NoStop}%
\bibitem [{\citenamefont {Schwabl}(2006)}]{Schwabl2006}%
  \BibitemOpen
  \bibfield  {author} {\bibinfo {author} {\bibfnamefont {F.}~\bibnamefont
  {Schwabl}},\ }\href {\doibase https://doi.org/10.1007/3-540-36217-7} {\emph
  {\bibinfo {title} {Statistical Mechanics}}}\ (\bibinfo  {publisher} {Springer
  Berlin Heidelberg},\ \bibinfo {year} {2006})\BibitemShut {NoStop}%
\bibitem [{\citenamefont {Zachos}\ \emph {et~al.}(2005)\citenamefont {Zachos},
  \citenamefont {Fairlie},\ and\ \citenamefont {Curtright}}]{Zachos2005}%
  \BibitemOpen
  \bibfield  {author} {\bibinfo {author} {\bibfnamefont {C.~K.}\ \bibnamefont
  {Zachos}}, \bibinfo {author} {\bibfnamefont {D.~B.}\ \bibnamefont {Fairlie}},
  \ and\ \bibinfo {author} {\bibfnamefont {T.~L.}\ \bibnamefont {Curtright}},\
  }\href {\doibase https://doi.org/10.1142/5287} {\emph {\bibinfo {title}
  {Quantum Mechanics in Phase Space}}}\ (\bibinfo  {publisher} {World
  Scientific Publishing Company},\ \bibinfo {year} {2005})\BibitemShut
  {NoStop}%
\bibitem [{\citenamefont {Tom{\'{e}}}\ and\ \citenamefont
  {de~Oliveira}(2015)}]{Tom2015}%
  \BibitemOpen
  \bibfield  {author} {\bibinfo {author} {\bibfnamefont {T.}~\bibnamefont
  {Tom{\'{e}}}}\ and\ \bibinfo {author} {\bibfnamefont {M.~J.}\ \bibnamefont
  {de~Oliveira}},\ }\href {\doibase https://doi.org/10.1007/978-3-319-11770-6}
  {\emph {\bibinfo {title} {Stochastic Dynamics and Irreversibility}}}\
  (\bibinfo  {publisher} {Springer International Publishing},\ \bibinfo {year}
  {2015})\BibitemShut {NoStop}%
\bibitem [{\citenamefont {Kang}\ and\ \citenamefont {Wilcox}(2017)}]{Kang2017}%
  \BibitemOpen
  \bibfield  {author} {\bibinfo {author} {\bibfnamefont {W.}~\bibnamefont
  {Kang}}\ and\ \bibinfo {author} {\bibfnamefont {L.~C.}\ \bibnamefont
  {Wilcox}},\ }\href {\doibase https://doi.org/10.1007/s10589-017-9910-0}
  {\bibfield  {journal} {\bibinfo  {journal} {Computational Optimization and
  Applications}\ }\textbf {\bibinfo {volume} {68}},\ \bibinfo {pages} {289}
  (\bibinfo {year} {2017})}\BibitemShut {NoStop}%
\bibitem [{\citenamefont {Griebel}(2006)}]{Griebel2005}%
  \BibitemOpen
  \bibfield  {author} {\bibinfo {author} {\bibfnamefont {M.}~\bibnamefont
  {Griebel}},\ }\href
  {http://wissrech.ins.uni-bonn.de/research/pub/griebel/focm.pdf} {\emph
  {\bibinfo {title} {Sparse grids and related approximation schemes for higher
  dimensional problems}}}\ (\bibinfo  {publisher} {Cambridge University
  Press},\ \bibinfo {year} {2006})\ pp.\ \bibinfo {pages}
  {106--161}\BibitemShut {NoStop}%
\bibitem [{\citenamefont {Moral}(2004)}]{DelMoral2004}%
  \BibitemOpen
  \bibfield  {author} {\bibinfo {author} {\bibfnamefont {P.~D.}\ \bibnamefont
  {Moral}},\ }\href {\doibase https://doi.org/10.1007/978-1-4684-9393-1} {\emph
  {\bibinfo {title} {Feynman-Kac Formulae}}}\ (\bibinfo  {publisher} {Springer
  New York},\ \bibinfo {year} {2004})\BibitemShut {NoStop}%
\bibitem [{\citenamefont {{\O}ksendal}(1998)}]{Oksendal1998}%
  \BibitemOpen
  \bibfield  {author} {\bibinfo {author} {\bibfnamefont {B.}~\bibnamefont
  {{\O}ksendal}},\ }\href {\doibase https://doi.org/10.1007/978-3-662-03620-4}
  {\emph {\bibinfo {title} {Stochastic Differential Equations}}}\ (\bibinfo
  {publisher} {Springer Berlin Heidelberg},\ \bibinfo {year}
  {1998})\BibitemShut {NoStop}%
\bibitem [{\citenamefont {Kraskov}\ \emph {et~al.}(2004)\citenamefont
  {Kraskov}, \citenamefont {St\"ogbauer},\ and\ \citenamefont
  {Grassberger}}]{Kraskov2004}%
  \BibitemOpen
  \bibfield  {author} {\bibinfo {author} {\bibfnamefont {A.}~\bibnamefont
  {Kraskov}}, \bibinfo {author} {\bibfnamefont {H.}~\bibnamefont
  {St\"ogbauer}}, \ and\ \bibinfo {author} {\bibfnamefont {P.}~\bibnamefont
  {Grassberger}},\ }\href {\doibase 10.1103/PhysRevE.69.066138} {\bibfield
  {journal} {\bibinfo  {journal} {Phys. Rev. E}\ }\textbf {\bibinfo {volume}
  {69}},\ \bibinfo {pages} {066138} (\bibinfo {year} {2004})}\BibitemShut
  {NoStop}%
\bibitem [{\citenamefont {Singh}\ and\ \citenamefont
  {Póczos}(2016)}]{Singh2016}%
  \BibitemOpen
  \bibfield  {author} {\bibinfo {author} {\bibfnamefont {S.}~\bibnamefont
  {Singh}}\ and\ \bibinfo {author} {\bibfnamefont {B.}~\bibnamefont
  {Póczos}},\ }\href@noop {} {\enquote {\bibinfo {title} {Finite-sample
  analysis of fixed-k nearest neighbor density functional estimators},}\ }
  (\bibinfo {year} {2016}),\ \Eprint {http://arxiv.org/abs/1606.01554}
  {arXiv:1606.01554 [math.ST]} \BibitemShut {NoStop}%
\bibitem [{\citenamefont {Ao}\ and\ \citenamefont {Li}(2022)}]{Ao2022}%
  \BibitemOpen
  \bibfield  {author} {\bibinfo {author} {\bibfnamefont {Z.}~\bibnamefont
  {Ao}}\ and\ \bibinfo {author} {\bibfnamefont {J.}~\bibnamefont {Li}},\ }\href
  {https://aaai-2022.virtualchair.net/poster_aaai8184} {\enquote {\bibinfo
  {title} {Entropy estimation via normalizing flow},}\ } (\bibinfo {year}
  {2022})\BibitemShut {NoStop}%
\bibitem [{\citenamefont {McMillan}(1965)}]{McMillan1965}%
  \BibitemOpen
  \bibfield  {author} {\bibinfo {author} {\bibfnamefont {W.~L.}\ \bibnamefont
  {McMillan}},\ }\href {\doibase 10.1103/PhysRev.138.A442} {\bibfield
  {journal} {\bibinfo  {journal} {Phys. Rev.}\ }\textbf {\bibinfo {volume}
  {138}},\ \bibinfo {pages} {A442} (\bibinfo {year} {1965})}\BibitemShut
  {NoStop}%
\bibitem [{\citenamefont {Ceperley}\ \emph {et~al.}(1977)\citenamefont
  {Ceperley}, \citenamefont {Chester},\ and\ \citenamefont
  {Kalos}}]{Ceperley1977}%
  \BibitemOpen
  \bibfield  {author} {\bibinfo {author} {\bibfnamefont {D.}~\bibnamefont
  {Ceperley}}, \bibinfo {author} {\bibfnamefont {G.~V.}\ \bibnamefont
  {Chester}}, \ and\ \bibinfo {author} {\bibfnamefont {M.~H.}\ \bibnamefont
  {Kalos}},\ }\href {\doibase 10.1103/PhysRevB.16.3081} {\bibfield  {journal}
  {\bibinfo  {journal} {Phys. Rev. B}\ }\textbf {\bibinfo {volume} {16}},\
  \bibinfo {pages} {3081} (\bibinfo {year} {1977})}\BibitemShut {NoStop}%
\bibitem [{\citenamefont {Carleo}\ \emph {et~al.}(2017)\citenamefont {Carleo},
  \citenamefont {Cevolani}, \citenamefont {Sanchez-Palencia},\ and\
  \citenamefont {Holzmann}}]{Carleo2017b}%
  \BibitemOpen
  \bibfield  {author} {\bibinfo {author} {\bibfnamefont {G.}~\bibnamefont
  {Carleo}}, \bibinfo {author} {\bibfnamefont {L.}~\bibnamefont {Cevolani}},
  \bibinfo {author} {\bibfnamefont {L.}~\bibnamefont {Sanchez-Palencia}}, \
  and\ \bibinfo {author} {\bibfnamefont {M.}~\bibnamefont {Holzmann}},\ }\href
  {\doibase 10.1103/physrevx.7.031026} {\bibfield  {journal} {\bibinfo
  {journal} {Physical Review X}\ }\textbf {\bibinfo {volume} {7}},\  (\bibinfo
  {year} {2017})}\BibitemShut {NoStop}%
\bibitem [{\citenamefont {Carleo}\ and\ \citenamefont
  {Troyer}(2017)}]{Carleo2017}%
  \BibitemOpen
  \bibfield  {author} {\bibinfo {author} {\bibfnamefont {G.}~\bibnamefont
  {Carleo}}\ and\ \bibinfo {author} {\bibfnamefont {M.}~\bibnamefont
  {Troyer}},\ }\href {\doibase 10.1126/science.aag2302} {\bibfield  {journal}
  {\bibinfo  {journal} {Science}\ }\textbf {\bibinfo {volume} {355}},\ \bibinfo
  {pages} {602} (\bibinfo {year} {2017})}\BibitemShut {NoStop}%
\bibitem [{\citenamefont {Dinh}\ \emph {et~al.}(2016)\citenamefont {Dinh},
  \citenamefont {Sohl-Dickstein},\ and\ \citenamefont {Bengio}}]{Dinh2016}%
  \BibitemOpen
  \bibfield  {author} {\bibinfo {author} {\bibfnamefont {L.}~\bibnamefont
  {Dinh}}, \bibinfo {author} {\bibfnamefont {J.}~\bibnamefont
  {Sohl-Dickstein}}, \ and\ \bibinfo {author} {\bibfnamefont {S.}~\bibnamefont
  {Bengio}},\ }\href {\doibase 10.48550/ARXIV.1605.08803} {\enquote {\bibinfo
  {title} {{Density estimation using Real NVP}},}\ } (\bibinfo {year}
  {2016})\BibitemShut {NoStop}%
\bibitem [{\citenamefont {Noé}\ \emph {et~al.}(2019)\citenamefont {Noé},
  \citenamefont {Olsson}, \citenamefont {Köhler},\ and\ \citenamefont
  {Wu}}]{Noe2019}%
  \BibitemOpen
  \bibfield  {author} {\bibinfo {author} {\bibfnamefont {F.}~\bibnamefont
  {Noé}}, \bibinfo {author} {\bibfnamefont {S.}~\bibnamefont {Olsson}},
  \bibinfo {author} {\bibfnamefont {J.}~\bibnamefont {Köhler}}, \ and\
  \bibinfo {author} {\bibfnamefont {H.}~\bibnamefont {Wu}},\ }\href {\doibase
  10.1126/science.aaw1147} {\bibfield  {journal} {\bibinfo  {journal}
  {Science}\ }\textbf {\bibinfo {volume} {365}},\ \bibinfo {pages} {eaaw1147}
  (\bibinfo {year} {2019})}\BibitemShut {NoStop}%
\bibitem [{\citenamefont {Papamakarios}\ \emph {et~al.}(2019)\citenamefont
  {Papamakarios}, \citenamefont {Nalisnick}, \citenamefont {Rezende},
  \citenamefont {Mohamed},\ and\ \citenamefont
  {Lakshminarayanan}}]{Papamakarios2019}%
  \BibitemOpen
  \bibfield  {author} {\bibinfo {author} {\bibfnamefont {G.}~\bibnamefont
  {Papamakarios}}, \bibinfo {author} {\bibfnamefont {E.}~\bibnamefont
  {Nalisnick}}, \bibinfo {author} {\bibfnamefont {D.~J.}\ \bibnamefont
  {Rezende}}, \bibinfo {author} {\bibfnamefont {S.}~\bibnamefont {Mohamed}}, \
  and\ \bibinfo {author} {\bibfnamefont {B.}~\bibnamefont {Lakshminarayanan}},\
  }\href@noop {} {\enquote {\bibinfo {title} {Normalizing flows for
  probabilistic modeling and inference},}\ } (\bibinfo {year} {2019}),\ \Eprint
  {http://arxiv.org/abs/1912.02762} {arXiv:1912.02762 [stat.ML]} \BibitemShut
  {NoStop}%
\bibitem [{\citenamefont {Grathwohl}\ \emph {et~al.}(2018)\citenamefont
  {Grathwohl}, \citenamefont {Chen}, \citenamefont {Bettencourt}, \citenamefont
  {Sutskever},\ and\ \citenamefont {Duvenaud}}]{Grathwohl2018}%
  \BibitemOpen
  \bibfield  {author} {\bibinfo {author} {\bibfnamefont {W.}~\bibnamefont
  {Grathwohl}}, \bibinfo {author} {\bibfnamefont {R.~T.~Q.}\ \bibnamefont
  {Chen}}, \bibinfo {author} {\bibfnamefont {J.}~\bibnamefont {Bettencourt}},
  \bibinfo {author} {\bibfnamefont {I.}~\bibnamefont {Sutskever}}, \ and\
  \bibinfo {author} {\bibfnamefont {D.}~\bibnamefont {Duvenaud}},\ }\href@noop
  {} {\enquote {\bibinfo {title} {Ffjord: Free-form continuous dynamics for
  scalable reversible generative models},}\ } (\bibinfo {year} {2018}),\
  \Eprint {http://arxiv.org/abs/1810.01367} {arXiv:1810.01367 [cs.LG]}
  \BibitemShut {NoStop}%
\bibitem [{\citenamefont {Dinh}\ \emph {et~al.}(2014)\citenamefont {Dinh},
  \citenamefont {Krueger},\ and\ \citenamefont {Bengio}}]{Dinh2014}%
  \BibitemOpen
  \bibfield  {author} {\bibinfo {author} {\bibfnamefont {L.}~\bibnamefont
  {Dinh}}, \bibinfo {author} {\bibfnamefont {D.}~\bibnamefont {Krueger}}, \
  and\ \bibinfo {author} {\bibfnamefont {Y.}~\bibnamefont {Bengio}},\
  }\href@noop {} {\enquote {\bibinfo {title} {Nice: Non-linear independent
  components estimation},}\ } (\bibinfo {year} {2014}),\ \Eprint
  {http://arxiv.org/abs/1410.8516} {arXiv:1410.8516 [cs.LG]} \BibitemShut
  {NoStop}%
\bibitem [{\citenamefont {Kingma}\ and\ \citenamefont
  {Dhariwal}(2018)}]{Kingma2018}%
  \BibitemOpen
  \bibfield  {author} {\bibinfo {author} {\bibfnamefont {D.~P.}\ \bibnamefont
  {Kingma}}\ and\ \bibinfo {author} {\bibfnamefont {P.}~\bibnamefont
  {Dhariwal}},\ }\href@noop {} {\enquote {\bibinfo {title} {Glow: Generative
  flow with invertible 1x1 convolutions},}\ } (\bibinfo {year} {2018}),\
  \Eprint {http://arxiv.org/abs/1807.03039} {arXiv:1807.03039 [stat.ML]}
  \BibitemShut {NoStop}%
\bibitem [{\citenamefont {Ardizzone}\ \emph {et~al.}(2018)\citenamefont
  {Ardizzone}, \citenamefont {Kruse}, \citenamefont {Wirkert}, \citenamefont
  {Rahner}, \citenamefont {Pellegrini}, \citenamefont {Klessen}, \citenamefont
  {Maier-Hein}, \citenamefont {Rother},\ and\ \citenamefont
  {Köthe}}]{Ardizzone2018}%
  \BibitemOpen
  \bibfield  {author} {\bibinfo {author} {\bibfnamefont {L.}~\bibnamefont
  {Ardizzone}}, \bibinfo {author} {\bibfnamefont {J.}~\bibnamefont {Kruse}},
  \bibinfo {author} {\bibfnamefont {S.}~\bibnamefont {Wirkert}}, \bibinfo
  {author} {\bibfnamefont {D.}~\bibnamefont {Rahner}}, \bibinfo {author}
  {\bibfnamefont {E.~W.}\ \bibnamefont {Pellegrini}}, \bibinfo {author}
  {\bibfnamefont {R.~S.}\ \bibnamefont {Klessen}}, \bibinfo {author}
  {\bibfnamefont {L.}~\bibnamefont {Maier-Hein}}, \bibinfo {author}
  {\bibfnamefont {C.}~\bibnamefont {Rother}}, \ and\ \bibinfo {author}
  {\bibfnamefont {U.}~\bibnamefont {Köthe}},\ }\href@noop {} {\enquote
  {\bibinfo {title} {Analyzing inverse problems with invertible neural
  networks},}\ } (\bibinfo {year} {2018}),\ \Eprint
  {http://arxiv.org/abs/1808.04730} {arXiv:1808.04730 [cs.LG]} \BibitemShut
  {NoStop}%
\bibitem [{\citenamefont {Albergo}\ \emph {et~al.}(2021)\citenamefont
  {Albergo}, \citenamefont {Kanwar}, \citenamefont {Racanière}, \citenamefont
  {Rezende}, \citenamefont {Urban}, \citenamefont {Boyda}, \citenamefont
  {Cranmer}, \citenamefont {Hackett},\ and\ \citenamefont
  {Shanahan}}]{Albergo2021}%
  \BibitemOpen
  \bibfield  {author} {\bibinfo {author} {\bibfnamefont {M.~S.}\ \bibnamefont
  {Albergo}}, \bibinfo {author} {\bibfnamefont {G.}~\bibnamefont {Kanwar}},
  \bibinfo {author} {\bibfnamefont {S.}~\bibnamefont {Racanière}}, \bibinfo
  {author} {\bibfnamefont {D.~J.}\ \bibnamefont {Rezende}}, \bibinfo {author}
  {\bibfnamefont {J.~M.}\ \bibnamefont {Urban}}, \bibinfo {author}
  {\bibfnamefont {D.}~\bibnamefont {Boyda}}, \bibinfo {author} {\bibfnamefont
  {K.}~\bibnamefont {Cranmer}}, \bibinfo {author} {\bibfnamefont {D.~C.}\
  \bibnamefont {Hackett}}, \ and\ \bibinfo {author} {\bibfnamefont {P.~E.}\
  \bibnamefont {Shanahan}},\ }\href {\doibase 10.1103/physrevd.104.114507}
  {\bibfield  {journal} {\bibinfo  {journal} {Physical Review D}\ }\textbf
  {\bibinfo {volume} {104}},\  (\bibinfo {year} {2021})}\BibitemShut {NoStop}%
\bibitem [{\citenamefont {Pawlowski}\ and\ \citenamefont
  {Urban}(2022)}]{Pawlowski2022}%
  \BibitemOpen
  \bibfield  {author} {\bibinfo {author} {\bibfnamefont {J.~M.}\ \bibnamefont
  {Pawlowski}}\ and\ \bibinfo {author} {\bibfnamefont {J.~M.}\ \bibnamefont
  {Urban}},\ }\href {\doibase 10.48550/ARXIV.2203.01243} {\enquote {\bibinfo
  {title} {Flow-based density of states for complex actions},}\ } (\bibinfo
  {year} {2022})\BibitemShut {NoStop}%
\bibitem [{\citenamefont {Kong}\ and\ \citenamefont
  {Chaudhuri}(2020)}]{Kong2020}%
  \BibitemOpen
  \bibfield  {author} {\bibinfo {author} {\bibfnamefont {Z.}~\bibnamefont
  {Kong}}\ and\ \bibinfo {author} {\bibfnamefont {K.}~\bibnamefont
  {Chaudhuri}},\ }\href@noop {} {\enquote {\bibinfo {title} {The expressive
  power of a class of normalizing flow models},}\ } (\bibinfo {year} {2020}),\
  \Eprint {http://arxiv.org/abs/2006.00392} {arXiv:2006.00392 [cs.LG]}
  \BibitemShut {NoStop}%
\bibitem [{\citenamefont {Teshima}\ \emph {et~al.}(2020)\citenamefont
  {Teshima}, \citenamefont {Ishikawa}, \citenamefont {Tojo}, \citenamefont
  {Oono}, \citenamefont {Ikeda},\ and\ \citenamefont {Sugiyama}}]{Teshima2020}%
  \BibitemOpen
  \bibfield  {author} {\bibinfo {author} {\bibfnamefont {T.}~\bibnamefont
  {Teshima}}, \bibinfo {author} {\bibfnamefont {I.}~\bibnamefont {Ishikawa}},
  \bibinfo {author} {\bibfnamefont {K.}~\bibnamefont {Tojo}}, \bibinfo {author}
  {\bibfnamefont {K.}~\bibnamefont {Oono}}, \bibinfo {author} {\bibfnamefont
  {M.}~\bibnamefont {Ikeda}}, \ and\ \bibinfo {author} {\bibfnamefont
  {M.}~\bibnamefont {Sugiyama}},\ }\href@noop {} {\enquote {\bibinfo {title}
  {Coupling-based invertible neural networks are universal diffeomorphism
  approximators},}\ } (\bibinfo {year} {2020}),\ \Eprint
  {http://arxiv.org/abs/2006.11469} {arXiv:2006.11469 [cs.LG]} \BibitemShut
  {NoStop}%
\bibitem [{\citenamefont {Reh}\ \emph {et~al.}(2021)\citenamefont {Reh},
  \citenamefont {Schmitt},\ and\ \citenamefont {G\"arttner}}]{Reh2021}%
  \BibitemOpen
  \bibfield  {author} {\bibinfo {author} {\bibfnamefont {M.}~\bibnamefont
  {Reh}}, \bibinfo {author} {\bibfnamefont {M.}~\bibnamefont {Schmitt}}, \ and\
  \bibinfo {author} {\bibfnamefont {M.}~\bibnamefont {G\"arttner}},\ }\href
  {\doibase 10.1103/PhysRevLett.127.230501} {\bibfield  {journal} {\bibinfo
  {journal} {Phys. Rev. Lett.}\ }\textbf {\bibinfo {volume} {127}},\ \bibinfo
  {pages} {230501} (\bibinfo {year} {2021})}\BibitemShut {NoStop}%
\bibitem [{\citenamefont {Brandstetter}\ \emph {et~al.}(2022)\citenamefont
  {Brandstetter}, \citenamefont {Worrall},\ and\ \citenamefont
  {Welling}}]{Brandstetter2022}%
  \BibitemOpen
  \bibfield  {author} {\bibinfo {author} {\bibfnamefont {J.}~\bibnamefont
  {Brandstetter}}, \bibinfo {author} {\bibfnamefont {D.}~\bibnamefont
  {Worrall}}, \ and\ \bibinfo {author} {\bibfnamefont {M.}~\bibnamefont
  {Welling}},\ }\href@noop {} {\enquote {\bibinfo {title} {{Message Passing
  Neural PDE Solvers}},}\ } (\bibinfo {year} {2022}),\ \Eprint
  {http://arxiv.org/abs/2202.03376} {arXiv:2202.03376 [cs.LG]} \BibitemShut
  {NoStop}%
\bibitem [{\citenamefont {Lu}\ \emph {et~al.}(2021)\citenamefont {Lu},
  \citenamefont {Jin}, \citenamefont {Pang}, \citenamefont {Zhang},\ and\
  \citenamefont {Karniadakis}}]{Lu2021}%
  \BibitemOpen
  \bibfield  {author} {\bibinfo {author} {\bibfnamefont {L.}~\bibnamefont
  {Lu}}, \bibinfo {author} {\bibfnamefont {P.}~\bibnamefont {Jin}}, \bibinfo
  {author} {\bibfnamefont {G.}~\bibnamefont {Pang}}, \bibinfo {author}
  {\bibfnamefont {Z.}~\bibnamefont {Zhang}}, \ and\ \bibinfo {author}
  {\bibfnamefont {G.~E.}\ \bibnamefont {Karniadakis}},\ }\href {\doibase
  10.1038/s42256-021-00302-5} {\bibfield  {journal} {\bibinfo  {journal}
  {Nature Machine Intelligence}\ }\textbf {\bibinfo {volume} {3}},\ \bibinfo
  {pages} {218–229} (\bibinfo {year} {2021})}\BibitemShut {NoStop}%
\bibitem [{\citenamefont {Li}\ \emph {et~al.}(2020)\citenamefont {Li},
  \citenamefont {Kovachki}, \citenamefont {Azizzadenesheli}, \citenamefont
  {Liu}, \citenamefont {Bhattacharya}, \citenamefont {Stuart},\ and\
  \citenamefont {Anandkumar}}]{Li2020}%
  \BibitemOpen
  \bibfield  {author} {\bibinfo {author} {\bibfnamefont {Z.}~\bibnamefont
  {Li}}, \bibinfo {author} {\bibfnamefont {N.}~\bibnamefont {Kovachki}},
  \bibinfo {author} {\bibfnamefont {K.}~\bibnamefont {Azizzadenesheli}},
  \bibinfo {author} {\bibfnamefont {B.}~\bibnamefont {Liu}}, \bibinfo {author}
  {\bibfnamefont {K.}~\bibnamefont {Bhattacharya}}, \bibinfo {author}
  {\bibfnamefont {A.}~\bibnamefont {Stuart}}, \ and\ \bibinfo {author}
  {\bibfnamefont {A.}~\bibnamefont {Anandkumar}},\ }\href@noop {} {\enquote
  {\bibinfo {title} {Fourier neural operator for parametric partial
  differential equations},}\ } (\bibinfo {year} {2020}),\ \Eprint
  {http://arxiv.org/abs/2010.08895} {arXiv:2010.08895 [cs.LG]} \BibitemShut
  {NoStop}%
\bibitem [{\citenamefont {Beck}\ \emph {et~al.}(2021)\citenamefont {Beck},
  \citenamefont {Becker}, \citenamefont {Grohs}, \citenamefont {Jaafari},\ and\
  \citenamefont {Jentzen}}]{Beck2021}%
  \BibitemOpen
  \bibfield  {author} {\bibinfo {author} {\bibfnamefont {C.}~\bibnamefont
  {Beck}}, \bibinfo {author} {\bibfnamefont {S.}~\bibnamefont {Becker}},
  \bibinfo {author} {\bibfnamefont {P.}~\bibnamefont {Grohs}}, \bibinfo
  {author} {\bibfnamefont {N.}~\bibnamefont {Jaafari}}, \ and\ \bibinfo
  {author} {\bibfnamefont {A.}~\bibnamefont {Jentzen}},\ }\href {\doibase
  10.1007/s10915-021-01590-0} {\bibfield  {journal} {\bibinfo  {journal}
  {Journal of Scientific Computing}\ }\textbf {\bibinfo {volume} {88}},\
  (\bibinfo {year} {2021})}\BibitemShut {NoStop}%
\bibitem [{\citenamefont {Raissi}\ \emph {et~al.}(2019)\citenamefont {Raissi},
  \citenamefont {Perdikaris},\ and\ \citenamefont {Karniadakis}}]{Raissi2019}%
  \BibitemOpen
  \bibfield  {author} {\bibinfo {author} {\bibfnamefont {M.}~\bibnamefont
  {Raissi}}, \bibinfo {author} {\bibfnamefont {P.}~\bibnamefont {Perdikaris}},
  \ and\ \bibinfo {author} {\bibfnamefont {G.}~\bibnamefont {Karniadakis}},\
  }\href {\doibase https://doi.org/10.1016/j.jcp.2018.10.045} {\bibfield
  {journal} {\bibinfo  {journal} {Journal of Computational Physics}\ }\textbf
  {\bibinfo {volume} {378}},\ \bibinfo {pages} {686} (\bibinfo {year}
  {2019})}\BibitemShut {NoStop}%
\bibitem [{\citenamefont {Feng}\ \emph {et~al.}(2022)\citenamefont {Feng},
  \citenamefont {Zeng},\ and\ \citenamefont {Zhou}}]{Feng2022}%
  \BibitemOpen
  \bibfield  {author} {\bibinfo {author} {\bibfnamefont {X.}~\bibnamefont
  {Feng}}, \bibinfo {author} {\bibfnamefont {L.}~\bibnamefont {Zeng}}, \ and\
  \bibinfo {author} {\bibfnamefont {T.}~\bibnamefont {Zhou}},\ }\href {\doibase
  10.2139/ssrn.4003881} {\bibfield  {journal} {\bibinfo  {journal} {SSRN
  Electronic Journal}\ } (\bibinfo {year} {2022}),\
  10.2139/ssrn.4003881}\BibitemShut {NoStop}%
\bibitem [{SM()}]{SM}%
  \BibitemOpen
  \href@noop {} {}\bibinfo {note} {{}See Supplemental Material for details on
  the TDVP derivation, estimation errors, algorithmic complexity, normalizing
  flows, the isotropic heat equation, phase-space dynamics and network
  hyperparameters.}\BibitemShut {Stop}%
\bibitem [{\citenamefont {Jaini}\ \emph {et~al.}(2019)\citenamefont {Jaini},
  \citenamefont {Selby},\ and\ \citenamefont {Yu}}]{Jaini2019}%
  \BibitemOpen
  \bibfield  {author} {\bibinfo {author} {\bibfnamefont {P.}~\bibnamefont
  {Jaini}}, \bibinfo {author} {\bibfnamefont {K.~A.}\ \bibnamefont {Selby}}, \
  and\ \bibinfo {author} {\bibfnamefont {Y.}~\bibnamefont {Yu}},\ }\href
  {\doibase 10.48550/ARXIV.1905.02325} {\enquote {\bibinfo {title}
  {Sum-of-squares polynomial flow},}\ } (\bibinfo {year} {2019}),\ \Eprint
  {http://arxiv.org/abs/1905.02325} {arXiv:1905.02325 [cs.LG]} \BibitemShut
  {NoStop}%
\bibitem [{\citenamefont {Schmitt}\ and\ \citenamefont
  {Heyl}(2020)}]{Schmitt2020}%
  \BibitemOpen
  \bibfield  {author} {\bibinfo {author} {\bibfnamefont {M.}~\bibnamefont
  {Schmitt}}\ and\ \bibinfo {author} {\bibfnamefont {M.}~\bibnamefont {Heyl}},\
  }\href {\doibase 10.1103/PhysRevLett.125.100503} {\bibfield  {journal}
  {\bibinfo  {journal} {Phys. Rev. Lett.}\ }\textbf {\bibinfo {volume} {125}},\
  \bibinfo {pages} {100503} (\bibinfo {year} {2020})}\BibitemShut {NoStop}%
\bibitem [{\citenamefont {Poirier}\ and\ \citenamefont
  {Geiger}(2016)}]{Poirier2016}%
  \BibitemOpen
  \bibfield  {author} {\bibinfo {author} {\bibfnamefont {D.~R.}\ \bibnamefont
  {Poirier}}\ and\ \bibinfo {author} {\bibfnamefont {G.~H.}\ \bibnamefont
  {Geiger}},\ }\href {\doibase https://doi.org/10.1007/978-3-319-48090-9_9}
  {\emph {\bibinfo {title} {Transport Phenomena in Materials Processing}}}\
  (\bibinfo  {publisher} {Springer International Publishing},\ \bibinfo {year}
  {2016})\BibitemShut {NoStop}%
\bibitem [{\citenamefont {Annaratone}(2010)}]{Annaratone2010}%
  \BibitemOpen
  \bibfield  {author} {\bibinfo {author} {\bibfnamefont {D.}~\bibnamefont
  {Annaratone}},\ }\href {\doibase https://doi.org/10.1007/978-3-642-03932-4}
  {\emph {\bibinfo {title} {Engineering Heat Transfer}}}\ (\bibinfo
  {publisher} {Springer Berlin Heidelberg},\ \bibinfo {year}
  {2010})\BibitemShut {NoStop}%
\bibitem [{\citenamefont {Black}\ and\ \citenamefont
  {Scholes}(1973)}]{Black1973}%
  \BibitemOpen
  \bibfield  {author} {\bibinfo {author} {\bibfnamefont {F.}~\bibnamefont
  {Black}}\ and\ \bibinfo {author} {\bibfnamefont {M.}~\bibnamefont
  {Scholes}},\ }\href {\doibase http://dx.doi.org/10.1086/260062} {\bibfield
  {journal} {\bibinfo  {journal} {Journal of Political Economy}\ }\textbf
  {\bibinfo {volume} {81}},\ \bibinfo {pages} {637} (\bibinfo {year}
  {1973})}\BibitemShut {NoStop}%
\bibitem [{\citenamefont {Guillaume}(2018)}]{Guillaume2018}%
  \BibitemOpen
  \bibfield  {author} {\bibinfo {author} {\bibfnamefont {T.}~\bibnamefont
  {Guillaume}},\ }\href {\doibase https://doi.org/10.1007/s10479-018-3001-1}
  {\bibfield  {journal} {\bibinfo  {journal} {Annals of Operations Research}\
  }\textbf {\bibinfo {volume} {281}},\ \bibinfo {pages} {229} (\bibinfo {year}
  {2018})}\BibitemShut {NoStop}%
\bibitem [{\citenamefont {Davies}(1989)}]{Davies1989}%
  \BibitemOpen
  \bibfield  {author} {\bibinfo {author} {\bibfnamefont {E.~B.}\ \bibnamefont
  {Davies}},\ }\href {\doibase 10.1017/cbo9780511566158} {\emph {\bibinfo
  {title} {Heat Kernels and Spectral Theory}}}\ (\bibinfo  {publisher}
  {Cambridge University Press},\ \bibinfo {year} {1989})\BibitemShut {NoStop}%
\bibitem [{\citenamefont {Presilla}\ \emph {et~al.}(1997)\citenamefont
  {Presilla}, \citenamefont {Onofrio},\ and\ \citenamefont
  {Patriarca}}]{Presilla1997}%
  \BibitemOpen
  \bibfield  {author} {\bibinfo {author} {\bibfnamefont {C.}~\bibnamefont
  {Presilla}}, \bibinfo {author} {\bibfnamefont {R.}~\bibnamefont {Onofrio}}, \
  and\ \bibinfo {author} {\bibfnamefont {M.}~\bibnamefont {Patriarca}},\ }\href
  {\doibase 10.1088/0305-4470/30/21/014} {\bibfield  {journal} {\bibinfo
  {journal} {Journal of Physics A: Mathematical and General}\ }\textbf
  {\bibinfo {volume} {30}},\ \bibinfo {pages} {7385} (\bibinfo {year}
  {1997})}\BibitemShut {NoStop}%
\bibitem [{\citenamefont {Lin}\ \emph {et~al.}(2021)\citenamefont {Lin},
  \citenamefont {Duan}, \citenamefont {Wang},\ and\ \citenamefont
  {Zhang}}]{Lin2021a}%
  \BibitemOpen
  \bibfield  {author} {\bibinfo {author} {\bibfnamefont {L.}~\bibnamefont
  {Lin}}, \bibinfo {author} {\bibfnamefont {J.}~\bibnamefont {Duan}}, \bibinfo
  {author} {\bibfnamefont {X.}~\bibnamefont {Wang}}, \ and\ \bibinfo {author}
  {\bibfnamefont {Y.}~\bibnamefont {Zhang}},\ }\href {\doibase
  10.1063/5.0048483} {\bibfield  {journal} {\bibinfo  {journal} {Chaos: An
  Interdisciplinary Journal of Nonlinear Science}\ }\textbf {\bibinfo {volume}
  {31}},\ \bibinfo {pages} {051105} (\bibinfo {year} {2021})}\BibitemShut
  {NoStop}%
\bibitem [{\citenamefont {Schmitt}\ and\ \citenamefont {Reh}(2022)}]{jvmc}%
  \BibitemOpen
  \bibfield  {author} {\bibinfo {author} {\bibfnamefont {M.}~\bibnamefont
  {Schmitt}}\ and\ \bibinfo {author} {\bibfnamefont {M.}~\bibnamefont {Reh}},\
  }\href {\doibase 10.21468/SciPostPhysCodeb.2} {\bibfield  {journal} {\bibinfo
   {journal} {SciPost Phys. Codebases}\ ,\ \bibinfo {pages} {2}} (\bibinfo
  {year} {2022})}\BibitemShut {NoStop}%
\bibitem [{\citenamefont {Heek}\ \emph {et~al.}(2020)\citenamefont {Heek},
  \citenamefont {Levskaya}, \citenamefont {Oliver}, \citenamefont {Ritter},
  \citenamefont {Rondepierre}, \citenamefont {Steiner},\ and\ \citenamefont
  {van {Z}ee}}]{flax2020github}%
  \BibitemOpen
  \bibfield  {author} {\bibinfo {author} {\bibfnamefont {J.}~\bibnamefont
  {Heek}}, \bibinfo {author} {\bibfnamefont {A.}~\bibnamefont {Levskaya}},
  \bibinfo {author} {\bibfnamefont {A.}~\bibnamefont {Oliver}}, \bibinfo
  {author} {\bibfnamefont {M.}~\bibnamefont {Ritter}}, \bibinfo {author}
  {\bibfnamefont {B.}~\bibnamefont {Rondepierre}}, \bibinfo {author}
  {\bibfnamefont {A.}~\bibnamefont {Steiner}}, \ and\ \bibinfo {author}
  {\bibfnamefont {M.}~\bibnamefont {van {Z}ee}},\ }\href
  {http://github.com/google/flax} {\enquote {\bibinfo {title} {{F}lax: A neural
  network library and ecosystem for {JAX}},}\ } (\bibinfo {year}
  {2020})\BibitemShut {NoStop}%
\bibitem [{\citenamefont {Bradbury}\ \emph {et~al.}(2018)\citenamefont
  {Bradbury}, \citenamefont {Frostig}, \citenamefont {Hawkins}, \citenamefont
  {Johnson}, \citenamefont {Leary}, \citenamefont {Maclaurin}, \citenamefont
  {Necula}, \citenamefont {Paszke}, \citenamefont {Vander{P}las}, \citenamefont
  {Wanderman-{M}ilne},\ and\ \citenamefont {Zhang}}]{Bradbury2018}%
  \BibitemOpen
  \bibfield  {author} {\bibinfo {author} {\bibfnamefont {J.}~\bibnamefont
  {Bradbury}}, \bibinfo {author} {\bibfnamefont {R.}~\bibnamefont {Frostig}},
  \bibinfo {author} {\bibfnamefont {P.}~\bibnamefont {Hawkins}}, \bibinfo
  {author} {\bibfnamefont {M.~J.}\ \bibnamefont {Johnson}}, \bibinfo {author}
  {\bibfnamefont {C.}~\bibnamefont {Leary}}, \bibinfo {author} {\bibfnamefont
  {D.}~\bibnamefont {Maclaurin}}, \bibinfo {author} {\bibfnamefont
  {G.}~\bibnamefont {Necula}}, \bibinfo {author} {\bibfnamefont
  {A.}~\bibnamefont {Paszke}}, \bibinfo {author} {\bibfnamefont
  {J.}~\bibnamefont {Vander{P}las}}, \bibinfo {author} {\bibfnamefont
  {S.}~\bibnamefont {Wanderman-{M}ilne}}, \ and\ \bibinfo {author}
  {\bibfnamefont {Q.}~\bibnamefont {Zhang}},\ }\href
  {http://github.com/google/jax} {\enquote {\bibinfo {title} {{JAX}: composable
  transformations of {P}ython+{N}um{P}y programs},}\ } (\bibinfo {year}
  {2018})\BibitemShut {NoStop}%
\bibitem [{\citenamefont {{J\"{u}lich Supercomputing Centre}}(2019)}]{JUWELS}%
  \BibitemOpen
  \bibfield  {author} {\bibinfo {author} {\bibnamefont {{J\"{u}lich
  Supercomputing Centre}}},\ }\href {http://dx.doi.org/10.17815/jlsrf-5-171}
  {\bibfield  {journal} {\bibinfo  {journal} {J. Large-Scale Res. Facilities}\
  }\textbf {\bibinfo {volume} {5}} (\bibinfo {year} {2019})}\BibitemShut
  {NoStop}%
\end{thebibliography}%
\fi

\iftrue
\clearpage
\newpage
\appendix

\begin{center}
    \textbf{Supplementary Materials}
\end{center}

\section{Derivation of the TDVP Equation}
The basic idea of a TDVP is to minimize the distance 
\begin{equation}
\mathcal{D}\left(p_{\theta(t)} + \dot{p}_{\theta(t)} \tau, p_{\theta(t)} + \sum_k \frac{\partial p_{\theta(t)}}{\partial \theta_k}\dot{\theta}_k \tau\right),
\end{equation}
between the evolved state at time $t+\tau$ and the updated network state with respect to the parameter updates $\dot{\theta}$ at each time $t$.
Here, we exemplarily derive 
Eq.~(5)
from the Hellinger distance $\mathcal{D}_H(p, q)=\int \sqrt{p(\mathbf{x})q(\mathbf{x})} d\mathbf{x}$, i.e.\ by maximizing the classical fidelity $F(p, q)=1-\mathcal{D}_H(p, q)$. As noted in the main text, an equivalent derivation is possible using the Kullback-Leibler divergence $\mathcal{D}_{KL}$ which leads to the same result.
For better readability, we drop the time index and continue with the optimality condition
\begin{equation}
\begin{aligned}
0 &= \frac{\partial}{\partial \dot{\theta}_k} F\left(p + \dot{p}\tau, p + \sum_{k^\prime}\frac{\partial p}{\partial \theta_{k^\prime}}\dot{\theta}_{k^\prime}\tau\right)\\
&=\frac{\partial}{\partial \dot{\theta}_k}\int  p(\textbf{x}) \sqrt{1 + a\tau + b \tau^2}d\mathbf{x},
\end{aligned}
\end{equation}
where $a$ and $b$ are given by
\begin{equation}
\begin{aligned}
a &= \frac{\partial \log p(\textbf{x})}{\partial t} + \sum_{k^\prime} \frac{\partial  \log p(\textbf{x})}{\partial \theta_{k^\prime}}\dot{\theta}_{k^\prime},\\
b &= \frac{\partial \log p(\textbf{x})}{\partial t}\sum_{k^\prime}\frac{\partial \log p(\textbf{x})}{\partial \theta_{k^\prime}}\dot{\theta}_{k^\prime}.
\end{aligned}
\end{equation}
Next we perform a second order expansion of the square root in the (small) time step $\tau$:
\begin{equation}
\sqrt{1 + a\tau + b\tau^2} = 1 + \frac{a\tau}{2} + \frac{\tau^2}{8} (4b - a^2) + \mathcal{O}(\tau^3).
\end{equation}
Using that the normalization of $p$ is conserved under the time evolution one finds that the term linear in $\tau$ vanishes:
\begin{equation}
\begin{aligned}
p(\textbf{x})a&=\int \left( \dot{p}(\textbf{x}) + \sum_{k^\prime} \frac{\partial p(\textbf{x})}{\partial \theta_{k^\prime}}\dot{\theta}_{k^\prime}\right)d\textbf{x}\\
&=\sum_{k^\prime} \dot{\theta}_{k^\prime}\frac{\partial}{\partial \theta_{k^\prime}}\int  p(\textbf{x})d\textbf{x}\\
&=\sum_{k^\prime}\dot{\theta}_{k^\prime}\frac{\partial}{\partial \theta_{k^\prime}}1\\
&=0.
\end{aligned}
\label{eqn:TDVP_linear_term_vanishes}
\end{equation}
Thus, the optimality condition becomes

\begin{equation}
\begin{aligned}
0 &=\frac{\partial}{\partial \dot{\theta}_k}\int \frac{p(\textbf{x})}{p(\textbf{x})^2} \left(4\dot{p}(\textbf{x})\sum_{k^\prime}\frac{\partial p(\textbf{x})}{\partial \theta_{k^\prime}}\dot{\theta}_{k^\prime}\right.\\
&\hphantom{=\frac{\partial}{\partial \dot{\theta}_k}\int d\textbf{x} \frac{p(\textbf{x})}{p(\textbf{x})^2}\left(\right.} \left. - (\dot{p}(\textbf{x}) + \sum_{k^\prime}\frac{\partial p(\textbf{x})}{\partial \theta_{k^\prime}}\dot{\theta}_{k^\prime})^2\right)d\textbf{x}\\
&=-\frac{\partial}{\partial \dot{\theta}_k}\int \frac{p(\textbf{x})}{p({\textbf{x}^2})}\left(\dot{p}(\textbf{x}) -\sum_{k^\prime} \frac{\partial p(\textbf{x})}{\partial \theta_{k^\prime}}\dot{\theta}_{k^\prime}\right)^2 d\textbf{x} \\
&=-\frac{\partial}{\partial \dot{\theta}_k}\int p(\textbf{x})\left(\frac{\partial \log p(\textbf{x})}{\partial t} - \sum_{k^\prime}\frac{\partial \log p(\textbf{x})}{\partial \theta_{k^\prime}}\dot{\theta}_{k^\prime}\right)^2 d\textbf{x}\\
&=2\int p(\textbf{x})\frac{\log p(\textbf{x})}{\partial \theta_k}\left(\frac{\partial \log p(\textbf{x})}{\partial t} - \sum_{k^\prime}\frac{\partial \log p(\textbf{x})}{\partial \theta_{k^\prime}}\dot{\theta}_{k^\prime}\right)d\textbf{x}.
\label{eqn:TDVP_derivation_penultimatestep}
\end{aligned}
\end{equation}
Dropping the factor of 2 we obtain an equation for the optimal parameter update $\dot{\theta}$:
\begin{equation}
\begin{aligned}
0 =& \underbrace{\int p(\textbf{x})\frac{\partial \log p(\textbf{x})}{\partial t} \frac{\partial \log p(\textbf{x})}{\partial \theta_k} d\textbf{x}}_{=F_k}\\
&-\sum_{k^\prime}\underbrace{\int p(\textbf{x}) \frac{\partial \log p(\textbf{x})}{\partial \theta_k}\frac{\partial \log p(\textbf{x})}{\partial \theta_{k^\prime}} d\textbf{x}}_{=S_{kk^\prime}}\dot{\theta}_{k^\prime}\\
=& \left\langle \frac{\partial \log p(\textbf{x})}{\partial t} \frac{\partial \log p(\textbf{x})}{\partial \theta_k}\right\rangle_{\mathbf{x}\sim p}\\
&-\sum_{k^\prime}\left\langle \frac{\partial \log p(\textbf{x})}{\partial \theta_k}\frac{\partial \log p(\textbf{x})}{\partial \theta_{k^\prime}}\right\rangle_{\mathbf{x}\sim p}\dot{\theta}_{k^\prime}.
\end{aligned}
\end{equation}
Importantly, we can now evaluate the integral by sampling according to the encoded probabilities $p(\textbf{x})$ since the integrand  is proportional to $p(\textbf{x})$ for both $F$ and $S$. This is a unique property of the distance measures $\mathcal{D}_H$ and $\mathcal{D}_{KL}$ while other distance measures, as for example the $L^2$ norm, do not lead to expressions of a form that can be efficiently evaluated from Monte Carlo samples.
The same derivation can be carried out without assuming normalization. In this case the form of $S$ and $F$ is altered to
\begin{equation}
\begin{aligned}
p(\textbf{x})&\rightarrow\frac{p(\textbf{x})}{\int p(\textbf{x}) d\textbf{x}}\\
\log p(\textbf{x})&\rightarrow\log p(\textbf{x}) - \log \int p(\textbf{x}) d{\textbf{x}}\\
\frac{\partial \log p(\textbf{x})}{\partial \theta_k}&\rightarrow\frac{\partial \log p(\textbf{x})}{\partial \theta_k}-\left\langle\frac{\partial \log p(\textbf{x})}{\partial \theta_k}\right\rangle_{\textbf{x}\sim p}\\
\frac{\partial \log p(\textbf{x})}{\partial t}&\rightarrow\frac{\partial \log p(\textbf{x})}{\partial t}-\left\langle\frac{\partial \log p(\textbf{x})}{\partial t}\right\rangle_{\textbf{x}\sim p} \,,
\end{aligned}
\end{equation}
where the last two lines are obtained using
\begin{equation}
\begin{aligned}
&\frac{\partial}{\partial \theta_k} \left( \log p(\textbf{x}) - \log \int p(\textbf{x}) d\textbf{x} \right)\\
=&\frac{\partial\log p(\textbf{x})}{\partial \theta_k} - \frac{\int \frac{\partial p(\textbf{x})}{\partial \theta_k} d\textbf{x}}{\int p(\textbf{x}') d\textbf{x}'}\\
=&\frac{\partial\log p(\textbf{x})}{\partial \theta_k} - \int \frac{p(\textbf{x})}{\int p(\textbf{x}')d\textbf{x}'} \frac{\partial \log p(\textbf{x})}{\partial \theta_k}d\textbf{x}\\
=&\frac{\partial\log p(\textbf{x})}{\partial \theta_k} - \left\langle\frac{\partial \log p(\textbf{x})}{\partial \theta_k}\right\rangle_{\textbf{x}\sim p}.
\end{aligned}
\end{equation}
Here, the log derivative trick was used in the third line. One may proceed similarly for the time derivative. Overall, this leaves us with a connected correlator structure instead of a simple correlator
\begin{equation}
\begin{aligned}
S_{kk'}=&\langle O_k(\textbf{x}) O_{k'}(\textbf{x})\rangle_{\mathbf{x}\sim p_\theta(t)}\\
- & \langle O_k(\textbf{x})\rangle_{\mathbf{x}\sim p_\theta(t)} \langle O_{k'}(\textbf{x})\rangle_{\mathbf{x}\sim p_\theta(t)},\\
F_k = & \langle O_k(\textbf{x}) \partial_t\log (p_{\theta(t)}(\textbf{x}))\rangle_{\mathbf{x}\sim p_{\theta(t)}}\\
- & \langle O_k(\textbf{x}) \rangle_{\mathbf{x}\sim p_{\theta(t)}} \langle\partial_t\log (p_{\theta(t)}(\textbf{x}))\rangle_{\mathbf{x}\sim p_{\theta(t)}}
\end{aligned}
\end{equation}
with $O_k$ the (logarithmic) variational derivative
\begin{equation}
O_k(\textbf{x}) = \partial_{\theta_k} \log (p_{\theta(t)}(\textbf{x})).
\end{equation}

We finally arrive at
\begin{equation}
\dot{\theta}_k = \tilde{S}^{-1}_{kk^\prime}F_{k^\prime}
\end{equation}
where the tilde is due to the fact that we cannot invert $S$ itself but rather need to regularize it because it is usually rank-deficient. One can show that the updates that were found are indeed maxima of the fidelity:

\begin{equation}
\begin{aligned}
&\frac{\partial^2}{\partial \dot{\theta}_k^2}F(p + \tau\dot{p}, p + \tau\sum_{k^\prime} \frac{\partial p}{\partial \theta_{k^\prime}}\dot{\theta}_{k^\prime})\\
=&\frac{\partial}{\partial \dot{\theta}_k}(F_k - \sum_{k^\prime}S_{kk^\prime}\dot{\theta}_{k^\prime})\\
=&-\sum_{k^\prime} S_{kk^\prime}\delta_{k^\prime k}\\
=&-S_{kk}\\
=& - \left\langle \left(O_k(\mathbf{x})-\left\langle O_k(\mathbf{x}) \right\rangle\right)^2 \right\rangle_{\textbf{x}\sim p}\\
\leq& 0.
\end{aligned}
\end{equation}

\section{Approximation Error}
The adjustment of the parameters to reflect changes in the probability carries an associated error, as the parameters can usually not be changed to perfectly reflect the time derivatives of all the sampled points used to estimate $F$.
The TDVP allows to quantify this error by estimating the residual
\begin{equation}
    \label{eqn:tdvp_error}
    r(t) = \frac{1}{N_s} \sum_i \left|\dot{p}(t, \mathbf{x}_i) - \sum_k \frac{\partial p(t, \mathbf{x}_i)}{\partial \theta_k} \dot{\theta}_k\right|^2.
\end{equation}

\section{Computational Complexity}
Here we compare the computational complexity of the explicit variational method that we are proposing to an iterative gradient descent based technique.

The operations carried out in the gradient descent based procedure are given in Alg.~\ref{alg:iterative_gradient_descent}.
To summarize, the algorithm computes time derivatives at the sampled points using the differential operator $\mathcal{F}$, which depends on the PDE under scrutiny (e.g. $\mathcal{F} = D\Delta_\mathbf{x}$ in the case of the heat equation, see
Eq.~(6))
. The time derivatives, weighted with some small time step $\tau$, are added to the current probability values $p_\theta(\mathbf{x}_i)$ and define the new regression targets. Using a loss function that is minimal when the encoded distribution agrees with the new regression targets, one searches for a new solution in the parameter space. Once a convergence criterion is met, the search stops and continues with the next time step.
This search can become costly, since the optimization problem is in general non-convex without convergence guarantees. It is therefore beneficial to avoid the iterative search using a closed form, as lined out in Alg.~\ref{alg:Explicit_scheme}.

\begin{algorithm}[H]
\caption{Iterative Gradient Descent}\label{alg:iterative_gradient_descent}
\begin{algorithmic}[1]
\Procedure{timeStep}{$p_\theta$}
    \State $K \gets \{\mathbf{x}_1,..,\mathbf{x}_N\}$ \Comment{obtain sample set $K$}
    \State $\partial_t p_\theta(\mathbf{x}_i)  \gets \mathcal{F}(p_\theta)(\mathbf{x}_i)$ \Comment{get time derivatives at each $\mathbf{x}$}
    \State
    \While{convergence criteria is not met}
    
    \State $L \gets \sum_i\mathcal{D}\left(p_\theta(\mathbf{x}_i), p_\theta(\mathbf{x}_i) + \tau \partial_t p_\theta(\mathbf{x}_i)\right)$ \Comment{Define Loss function, Derivative acts on first argument in $\mathcal{D}$}
    \State $\theta \gets \theta + \eta \nabla_\theta L$ \Comment{Gradient descent step}
    \EndWhile
    \State    
    \State \textbf{return} $\theta$
\EndProcedure
\end{algorithmic}
\end{algorithm}

\begin{algorithm}[H]
\caption{Explicit second order scheme}\label{alg:Explicit_scheme}
\begin{algorithmic}[1]
\Procedure{timeStep}{$p_\theta$}
    \State $K \gets \{\mathbf{x}_1,..,\mathbf{x}_N\}$ \Comment{obtain sample set $K$}
    \State $\partial_t p_\theta(\mathbf{x}_i)  \gets \mathcal{F}(p_\theta)(\mathbf{x}_i)$ \Comment{get time derivatives at each $\mathbf{x}$}
    \State
    \State $S_{kk'} \gets 1/N  \sum_{i\in S} O_k(\textbf{x}_i) O_{k'}(\textbf{x}_i)$
    \State $F_k \gets 1/N \sum_{i\in S} O_k(\textbf{x}_i) \partial_t \log\left(p_\theta(\mathbf{x}_i)\right)$
    \State $\dot{\theta}_k \gets S_{kk'}^{-1}F_{k'}$
    \State $\theta \gets \theta + \tau \dot{\theta}$
    
    \State
    \State \textbf{return} $\theta$
\EndProcedure
\end{algorithmic}
\end{algorithm}

The proposed algorithm also relies on obtaining samples from the distribution for which then time derivatives are computed. Beyond that, we require the logarithmic variational derivatives of the probabilities $O_k(\mathbf{x}) = \partial_{\theta_k} \log\left(p_\theta(\mathbf{x})\right)$. Obtaining these derivatives is of similar computational cost compared to a single gradient descent step in Alg.~\ref{alg:iterative_gradient_descent}, as such a step also requires the differentiation of the probabilities at all sample positions. However, using the explicit scheme, this operation needs to be carried out only once. Additionally, we need to add these derivatives together, as shown in lines 5 and 6 of Alg.~\ref{alg:Explicit_scheme}, an operation with negligible computational cost. To arrive at the time derivatives of the parameters $\dot{\theta}$, we need to invert the $S$ matrix, which has cubic computational cost in the number of network parameters. In practice, there are many more computationally efficient ways than actually computing the inverse, which allow to reduce the computational cost of this step. While this step limits the number of network parameters that can practically be used, and thus the expressivity of the ansatz, we have not found this to be a limiting factor in the application considered in this work.

The runtimes of the examples we presented all lie below half an hour on a single NVIDIA A100 GPU.


\section{Normalizing Flows}

\begin{figure}[t]
    \centering
    \includegraphics[width=.95\linewidth]{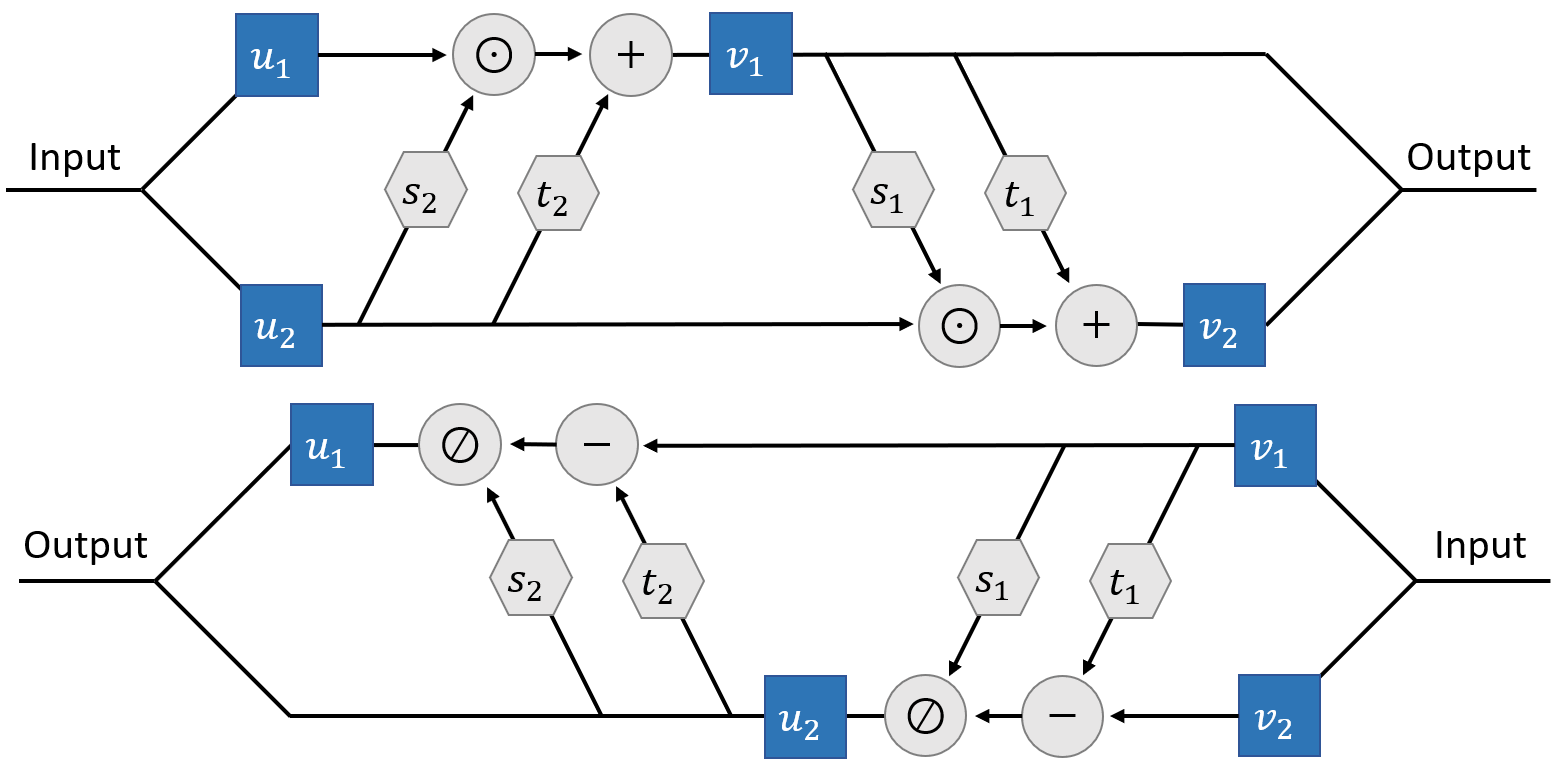}
    \caption{Illustration of a coupling block. Top: Forward pass, corresponding to Eq.~\eqref{eqn:coupling_forward}. Bottom: Inverse pass, corresponding to Eq.~\eqref{eqn:coupling_inverse}. 
    }
    \label{fig:couplingBlock}
\end{figure}
For the definition of the normalizing flow we use the Real-NVP type coupling blocks introduced in 
\cite{Dinh2016} 
and depicted in Fig.~\ref{fig:couplingBlock}. Each coupling block $\boldsymbol{\varphi}_i$ splits the input into two parts $u_1$ and $u_2$ which is done in a random but fixed way. The transformations in a single coupling block are defined by four networks $s_1$, $t_1$, $s_2$, $t_2$ that change the input as follows:
\begin{equation}
\begin{split}
    \label{eqn:coupling_forward}
    \mathbf{v}_1 & = \mathbf{u}_1 \odot \exp(\mathbf{s}_2(\mathbf{u}_2)) + \mathbf{t}_2(\mathbf{u}_2),\\
    \mathbf{v}_2 & = \mathbf{u}_2 \odot \exp(\mathbf{s}_1(\mathbf{v}_1)) + \mathbf{t}_1(\mathbf{v}_1).
\end{split}
\end{equation}
The inverse of this transformation is given by:
\begin{equation}
\begin{split}
    \label{eqn:coupling_inverse}
    \mathbf{u}_2 & = (\mathbf{v}_2 - \mathbf{t}_1(\mathbf{v}_1)) \odot \exp(-\mathbf{s}_1(\mathbf{v}_1)),\\
    \mathbf{u}_1 & = (\mathbf{v}_1 - \mathbf{t}_2(\mathbf{u}_2)) \odot \exp(-\mathbf{s}_2(\mathbf{u}_2)).
\end{split}
\end{equation}

Here $\odot$ means element-wise multiplication. The networks $\mathbf{s}$ and $\mathbf{t}$ are built equivalently as two layer feed-forward networks with half as many nodes in each layer as there are dimensions. In some cases we found it useful to not include the additive $\mathbf{t}$ networks. Additionally, we allowed the network to adjust the mean $\boldsymbol{\mu}$ and the covariance matrix $\Sigma$ of the distribution in latent space directly, potentially along with parameters $\vartheta$ of the distribution, e.g. $\nu$ in the case of the Student-$t$, such that
\begin{equation}
    \pi = \pi(\boldsymbol{\mu}_\theta, \Sigma_\theta, \vartheta_\theta).    
\end{equation}
We parameterize $\Sigma$ using either the Cholesky decomposition or by setting $\Sigma=\mathds{1}+AA^T$, where we found the latter to be more stable numerically for simulating the heat equation.
Network details are listed in in Table~\ref{tab:network_details}.
\onecolumngrid

\begin{table}[h]
    \centering
    \begin{tabular}{|c|c|c|c|c|c|c|c|}
    
    \hline
    Figure & Input Dim. & \# Coupling Blocks & \# Layers & Net $t$ & \# Parameters & \# Samples & $\pi$ \\ \hline
    Fig. 2 & 8 & 4 & 2 & No & 392 & 10.000 & Student-$t$ / Gauss \\
    Fig. 3 & 6 & 4 & 2 & Yes & 234 & 10.000 & Gauss \\
    \hline
    \end{tabular}
    \caption{\label{tab:network_details} Hyperparameters that were used for the different figures in the main text.}
\end{table}
\twocolumngrid

\begin{figure}
    \centering
    \includegraphics[width=.95\linewidth]{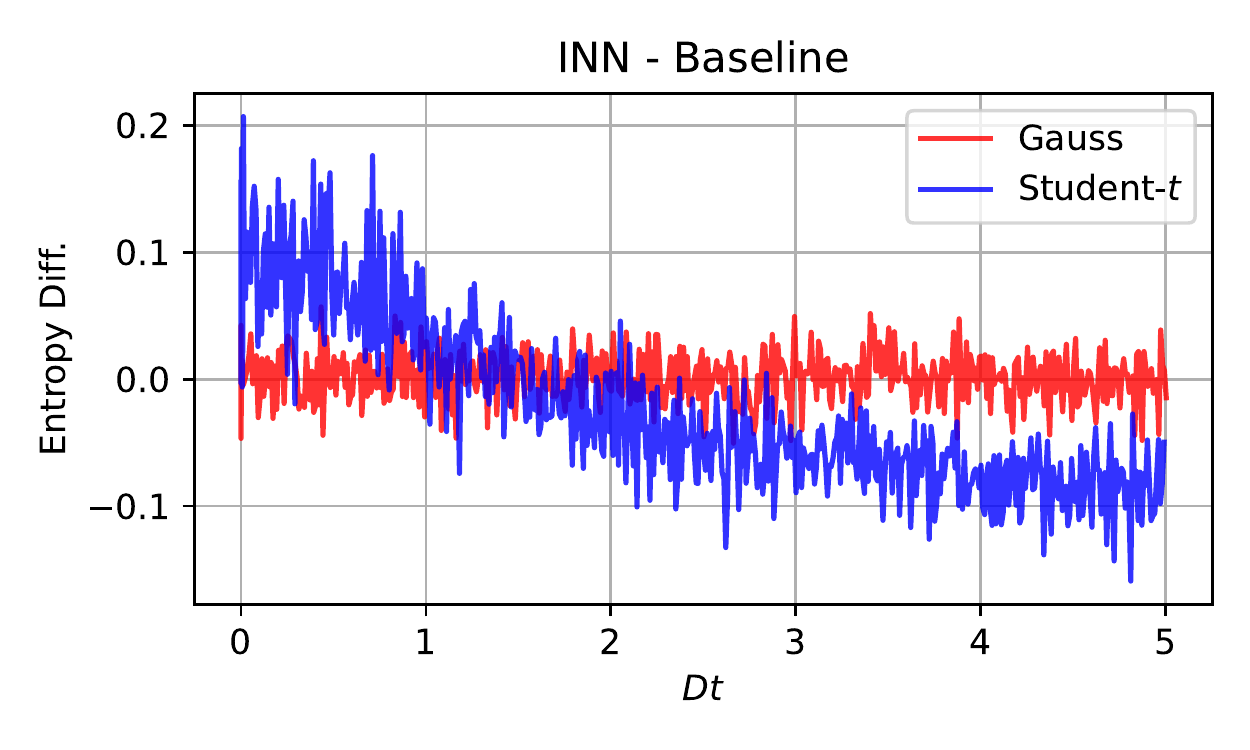}
    \caption{Differences between entropy estimates of the INN compared to the baseline methods from 
    Fig.~2
    in the main text.
    Systematic deviations are visible for the case of the Student-$t$, which we attribute to problems with the baseline method as layed out below.}
    \label{fig:diffusion_diff}
\end{figure}
\section{Isotropic Heat Equation as a 1D Problem}
Here we describe the procedure with which the reference data for 
Fig.~2 in the main text 
was obtained in the case of a Student-$t$ initial distribution and give an explanation for the slight discrepancies observed in 
Fig.~2
in the main text.

The heat equation
\begin{equation}
    \partial_t p(t, \mathbf{x}) = D \Delta p(t, \mathbf{x})
\end{equation}
can be recast as a 1D problem if the initial condition $p(t, \mathbf{x}) = u(\mathbf{x})$ features a spherical symmetry. This is the case if it is fully described by a mean $\boldsymbol{\mu}$ and a covariance matrix $\Sigma$, as this allows to rescale coordinates such that the new distribution obeys $\boldsymbol{\mu}=0$ and $\Sigma=\mathds{1}$ enabling us to write $p(t, \mathbf{x})$ as $p(t, r)$, where $r=|\mathbf{x}|$.

Then, the spherical form of the Laplacian may be exploited
\begin{equation}
\label{eqn:spherical_Laplace}
    \Delta = \partial_r^2 + \frac{d - 1}{r} \partial_r,
\end{equation}
where $d$ is the dimension of the distribution. The evolution of the distribution can then be solved using finite differences on a 1D grid. Note however, that there are caveats associated with this procedure as the second term of Eq.~\eqref{eqn:spherical_Laplace} has a divergence for $r\rightarrow0$. For the diffusion cases we considered, the distribution $p(t, r)$ has a maximum at $r=0$, irrespective of the time $t$, implying that also the numerator $\partial_t p(t, r)$ vanishes. This necessitates the use of L'H\^{o}spital's rule to write
\begin{equation}
    \lim_{r\rightarrow 0} \frac{\partial_r p}{r} = \frac{\partial_r^2 p}{\partial_r r}=\partial_r^2 p.
\end{equation}
We work with equidistant grid cells of size $\delta = 4\cdot 10^{-3}$ and set a cutoff at $r=100$. We employ L'H\^{o}spital approximation for the first 10 grid cells, i.e. for $r\in [0, 10\delta]$, which we found to be necessary for numerical stability. This implies however, that also the reference data is not free of approximations, which may be particularly interesting as we did not observe the INN curve to come closer to the reference data when increasing the network size. This could be viewed as an indication that it is not necessarily the INN whose curve is deviating from the true entropy, but rather the data obtained from the 1D grid-method described in this section. The difference between INN and grid-based result is shown in Fig.~\ref{fig:diffusion_diff}.

\section{Phase Space Evolution}
The phase space evolution of the example discussed in 
Fig.~3
of the main text is governed by the Hamiltonian $H$, which we choose to be
\begin{equation}
H = \sum_i\frac{1}{2}(m\omega^2 x_i^2 + p_i^2/m)+ k\sum_i (x_i - x_{(i+1) \% N})^2,
\end{equation}
such that $k$ gives the strength of the coupling between oscillators.
The resulting phase space flow is given by
\begin{equation}
    \begin{split}
    \dot{x}_i =& \partial_{p_i} H,\\
    \dot{p}_i =& -\partial_{x_i} H.\\
    \end{split}
\end{equation}
If one considers damping in phase space, the following Fokker-Planck equation is obtained
\cite{Presilla1997}
\begin{equation}
\begin{split}
\partial_t \rho(t, \mathbf{x}, \mathbf{p}) = &\left[-\partial_\mathbf{p}H \cdot \partial_\mathbf{x} + \partial_\mathbf{x}H\cdot\partial_\mathbf{p}+ \right.\\
&\left. \gamma\left(\mathbf{p}\cdot\partial_\mathbf{p} + m k_B \mathsmaller{\sum}_i T_i\partial_{p_i}^2 \right)\right] \rho(t, \mathbf{x}, \mathbf{p}).
\end{split}
\end{equation}

\fi

\end{document}